\documentclass[12pt]{article}

\usepackage{amsthm, amsmath, amssymb, tikz, ifthen, lscape}
\usepackage[scale=0.765]{geometry}
\usepackage{doi}
\usepackage{suetterl}
\usepackage[T1]{fontenc}

\DeclareMathOperator{\ZFC}{ZFC}
\DeclareMathOperator{\MA}{MA}
\DeclareMathOperator{\cf}{cf}

\DeclareMathOperator{\dom}{dom}
\DeclareMathOperator{\ran}{ran}
\DeclareMathOperator{\otyp}{otyp}

\DeclareMathOperator{\ld}{ld}

\newcommand{\con}{{Contradiction!}}
\newcommand{\arrows}{\rightarrow}
\newcommand{\doesntarrow}{\not\rightarrow}
\newcommand{\surj}{\; \rightarrow \mspace{-26mu} \longrightarrow}
\newcommand{\soast}[2]{\{#1|#2\}}
\newcommand{\soafft}[2]{^{#1} \mspace{-4mu} #2}
\newcommand{\card}[1]{{\scriptscriptstyle \overline{\overline{\textstyle #1}}}}
\newcommand{\pwi}[2]{#1 \mspace{-1mu} ^{``} \mspace{-2mu} #2}
\newcommand{\rb}[1]{\raisebox{1.5ex}[-1.5ex]{#1}}

\newcommand{\?}{\mspace{-1mu}}

\newcommand{\yq}{\mspace{9mu}}
\newcommand{\qc}{\mspace{1mu}}

\newcommand{\zbeat}{{\,\mapsto_r}}
\newcommand{\obeat}{{\,\mapsto_b}}
\newcommand{\tbeat}{{\,\mapsto_g}}
\newcommand{\beat}{\mapsto}

\newcommand{\uzbeat}{{\,-_r}}
\newcommand{\uobeat}{{\,-_b}}
\newcommand{\utbeat}{{\,-_g}}
\newcommand{\zinh}{N^-_R}
\newcommand{\zonh}{N^+_R}
\newcommand{\oinh}{N^-_B}
\newcommand{\oonh}{N^+_B}
\newcommand{\tinh}{N^-_G}
\newcommand{\tonh}{N^+_G}
\newcommand{\zb}{red}
\newcommand{\ob}{blue}

\newcommand{\ind}{I}
\newcommand{\On}{\Omega}
\newcommand{\this}{{\footnotesize This paper!}}\newtheorem{notation}{Notation}[section]
\newcommand{\ox}{}

\newtheorem{definition}[notation]{Definition}
\newtheorem{claim}[notation]{Claim}
\newtheorem{fact}[notation]{Fact}
\newtheorem{theorem}[notation]{Theorem}
\newtheorem{corollary}[notation]{Corollary}
\newtheorem{lemma}[notation]{Lemma}
\newtheorem{question}[notation]{Question}
\newtheorem{conjecture}[notation]{Conjecture}
\newtheorem{problem}[notation]{Problem}

\renewenvironment{proof}{Proof.}{\hspace*{\fill}$\blacksquare$ \vspace{2cc}} 
\newenvironment{nestedproof}{Proof of Claim:}{\hspace*{\fill}$\blacksquare\text{(Claim)}$ \vspace{2cc}}

\definecolor{ZERO}{rgb}{1,0,0}
\definecolor{ONE}{rgb}{0,0,1}
\definecolor{TWO}{rgb}{0,1,0}

\newcommand{\zero}{red}
\newcommand{\one}{blue}
\newcommand{\two}{green}

\definecolor{grey}{rgb}{.47,.47,.47}
\definecolor{very light grey}{rgb}{.93,.93,.93}
\definecolor{dark grey}{rgb}{.2,.2,.2}

\newcounter{tiefe}

\newenvironment{aufz}{\begin{list}{\ifthenelse{\value{tiefe}<3}{\ifthenelse{\value{tiefe}=1}{$\diamondsuit$}{$\bullet$}}
{\ifthenelse{\value{tiefe}=3}{$\diamond$}{\cdot}}}{\leftmargin=5dd}{\addtocounter{tiefe}{1}}}{\addtocounter{tiefe}{-1}\end{list}}

\usetikzlibrary{fadings}
\usetikzlibrary{patterns}

\begin{tikzfadingfrompicture}[name=titlefade] \node [text=transparent!20]
{\LARGE
Idiosynchromatic Poetry};
\end{tikzfadingfrompicture}
\date{February 8, 2015}
\title{\begin{tikzpicture} \fill [black!0] (-7,-1) rectangle (7,1); ]
(-7,-1) rectangle (7,1);
\shade[path fading=titlefade,fit fading=false, left color=red, right color=blue, middle color=green]
(-7,-1) rectangle (7,1);
\end{tikzpicture}}

\author{Thilo Weinert}

\begin{document}

\maketitle

\begin{abstract}

We prove equivalences of asymmetric partition relations involving natural numbers and products of weakly compact cardinals $\kappa$, infinite cardinals $\lambda<\kappa$ and natural numbers to certain classes of finitary problems in the theory of edge-coloured digraphs. We are able to determine three classes of Ramsey numbers exactly, typical examples in these classes are $r(\omega^2 2, 3)$, $r(\kappa\lambda 2, 3)$\ox\ and $r(\kappa\lambda 3, 3)$. We moreover provide general upper bounds for $r(\omega^2 m, 3)$ and $r(\kappa\lambda m, n)$.


{\bfseries Mathematical Subject Classification of the AMS:} 03E02, 05C20, 05C63, 05D10.
\end{abstract}

\section{Introduction}

\footnotetext{This title is a tribute both to my PhD advisor who regarded this work at an early stage as poetry---the author concedes but hopes it to be good poetry---and to David Schrittesser who came up with the word ``idiosynchromatic".}
\footnotetext{This is the coloured postprint version, the original greyscale version appeared as \cite{014W0}.}

The partition calculus as introduced by Erd\H{o}s and Rado in \cite{956ER0} is of interest both to set theorists and combinatorialists. While the former are usually more interested in statements with consistency strength which are hence independent over $\ZFC$ the latter are normally interested in finitary problems. In between falls a class of problems the \ae sthetic appeal of which is in the author's opinion usually underestimated---partition relations between countable ordinals.

\begin{notation}

In this paper we use the terminology which is nowadays pretty common. Ordinals are von-Neumann-ordinals so in particular an ordinal is identical to the set of its ordinal predecessors. Therefore $\alpha \setminus \beta$ denotes the half-open interval from $\beta$ to $\alpha$ containing the former but not the latter. $\On$ denotes the class of all ordinals. $[\alpha]^\beta$ denotes the set of subsets of $\alpha$ with order type $\beta$. $\pwi{c}{X}$ for $X \subset \dom(c)$ is defined to be $\soast{y \in \ran(c)}{\exists x \in X : c(x) = y}$ i.e. the pointwise image of $X$ under the function $c$. Recall that the asymmetric partition relation $\alpha \arrows (\beta, n)$ means $\forall c : [\alpha]^2 \longrightarrow 2\exists X \in [\alpha]^\beta : \pwi{c}{X} = \{0\} \vee \exists X \in [\alpha]^n : \pwi{c}{X} = \{1\}$, i.e. for every colouring of the unordered pairs of ordinals less than $\alpha$ with two colours we find a subset of $\alpha$ of order type $\beta$ homogeneous in colour $0$ or a subset of $\alpha$ of order type $n$ homogeneous in colour $1$.

Moreover we employ the common notation for Ramsey numbers also in this more general case. So $r(\beta, n) = \alpha$ means that $\alpha$ is the least ordinal such that $\alpha \arrows (\beta, n)$ holds true.

When $v, w$ are vertices of a directed graph $v\zbeat w$ means that there is a \zero\ arrow from $v$ to $w$ while $v \obeat w$ says that there is a \one\ arrow from $v$ to $w$, $v \tbeat w$ claims the existence of a \two\ arrow from $v$ to $w$\ox\ and finally $v \beat w$ means that there is an arrow of some kind from $v$ to $w$. Sometimes we want to talk about there being an arrow of a certain colour between $v$ and $w$ without claiming anything about its direction. For this we write $v\uzbeat w$---there is a \zero\ arrow between $v$ and $w$, $v\uobeat w$---there is a \one\ arrow between $v$ and $w$\ox\ or $v\utbeat w$---there is a \two\ arrow between $v$ and $w$. The \zero\ in-neighborhoods of a vertex $v$ will be denoted by $\zinh(v)$, the \zero\ out-neighborhoods by $\zonh(v)$\ox\ and $\oinh(v)$, $\oonh(v)$, $\tinh(v)$\ox\ and $\tonh(v)$ have their expected denotations. The set of vertices independent from $v$ will be denoted by $\ind(v)$.

In this paper, a digraph has at most one arrow between any two vertices and never any from a vertex to itself.

For the graph-theoretic terminology we follow \cite{009BG0}. Some of the notation was inspired by \cite{009GS0}.

\end{notation}

People working in this area seem to have been mainly interested in evaluating partition relations of the form $\alpha \arrows (\alpha, n)$ for countable $\alpha$ and natural $n$. There are three cases in which partition relations of the form $\alpha \arrows (\beta, n)$ for $\beta < \alpha$ both countable and $n$ natural have been analysed in some generality. This was the case for $\alpha, \beta$ finite powers of $\omega$, for $\alpha, \beta$ finite multiples of $\omega$ and for $\alpha, \beta$ natural---the last case of course being by far the most popular. In this last case people went to great lengths to calculate the $r(m, n)$'s, see for example \cite{995MR0}! There are a few results which do not fall into one of these categories, for example in \cite{969M0}, \cite{969HS0}, \cite{969HS1}\ox\ and \cite{969HS2}.

\begin{definition}

A \emph{coloured digraph} is a triple $\langle V, A, c \rangle$ such that $V$ is a set( of vertices), $A$ is a binary irreflexive asymmetric relation on $V$---a set of arcs---and $c : A \longrightarrow \ran(c)$ is a colouring of A.

\end{definition}

We are in the following going to be concerned with arc-coloured digraphs, i.e. directed graphs with coloured arcs. We will be especially interested in three-coloured digraphs, i.e. in the case where $\ran(c)=3$. Remember that a complete digraph is called a \emph{tournament}. For brevity in the following we will speak of \emph{triples} instead of three-person-three-coloured-tournaments.

In \cite{956ER0} Erd\H{o}s and Rado prove the following theorem\footnote{The author changed the notation partially to his own. This Theorem was first published as Theorem 25 in ``A Partition Calculus in Set Theory'' in the Bulletin of the American Mathematical Society in No. 62, pp. 427--489, published by the American Mathematical Society, \textcopyright P\'{a}l Erd\H{o}s and Richard Rado, 1956.}

\begin{theorem}

Let $m, n \in \omega \setminus 2$ and denote by $l_0 = l_0(m, n)$ the least finite number $l$ possessing the following property.

Property $P_{mn}$. Whenever $\rho(\lambda, \mu) < 2$ for $\{\lambda, \mu\} \in [l]^2$, then there is either $\{\lambda_0, \dots, \lambda_{m - 1}\} \in [l]^m$ such that $\rho(\lambda_\alpha, \lambda_\beta) = 0$ for $\alpha < \beta < m$ or there is $\{\lambda_0, \dots, \lambda_{n - 1}\} \in [l]^n$ such that $\rho(\lambda_\alpha, \lambda_\beta) = 1$ for $\{\alpha, \beta\} \in [n]^2$. Then $\omega l_0 \arrows (m, \omega n)^2$ but $\gamma \doesntarrow (m, \omega n)^2$ for $\gamma < \omega l$. Moreover, if $l_1 \arrows (m, m, n)^2$, then $l_0 \leqslant l_1$.

\label{erdosrado}

\end{theorem}

Note, however, that a possible alternative formulation of all but the last sentence of this theorem would nowadays be the following.

\begin{theorem}
The partition relation $\omega l \arrows (\omega m, n)$ holds true if and only if for every digraph $C$ on $l$ vertices contains an independent set of size $m$ or there is a subtournament $S$ of $C$ induced by a set of $n$ vertices such that all triples in $S$ are transitive.
\label{erdosradomodern}
\end{theorem}

These two theorems are equivalent by a well-known fact.

\begin{fact}
If a tournament has a cycle then it has one of length $3$.
\end{fact}

Later theorem \ref{erdosrado} was generalized by Baumgartner, see \cite{974B0}, in the following way:

\begin{theorem}
Let $\kappa$ be any infinite cardinal. The partition relation $\kappa l \arrows (\kappa m, n)$ holds true if and only if for every digraph $C$ on $l$ vertices contains an independent set of size $m$ or there is a subtournament $S$ of $C$ induced by a set of $n$ vertices such that all triples in $S$ are transitive.
\end{theorem}

\subsection{Added in Proof}
The author believed for quite some time that the discovery that $r(\omega^2 2, 3) = \omega^2 10$ was due to him. Only later he discovered that it already appeared in \cite{969HS2}. No proof was given though and there was also neither a formulation of a finitary problem equivalent to the calculation of the $r(I_m, A_n)$'s nor an upper bound for the $r(I_m, A_3)$'s given.

\section{Determining a Ramsey number}
In this section we mainly build on work of P\'al Erd\H{o}s and Richard Rado in the aforementioned \cite{956ER0} and Ernst Specker in \cite{957S0}.
The following eclectic definition is justified by theorem \ref{main} the proof of which explains the background of it.

\begin{definition}
We identify $0$ with \zero, $1$ with \one\ and $2$ with \two\ arrows.\footnote{This is a coloured postprint version.}
A triple is called \emph{agreeable} if and only if it is one of the following.

\label{listdefinition}
\begin{tikzpicture}[
zero/.style={ultra thick, solid, ZERO},
one/.style={ultra thick, solid, ONE},
two/.style={ultra thick, solid, TWO}
]

\pgftransformscale{0.023} 

\path
node (up) at (0, 71) {}
node (left) at (-41,0) {}
node (right) at (41,0) {};
\draw [->] (up) [zero] -- (left);
\draw [->] (left) [zero] -- (right);
\draw [->] (up) [zero] -- (right);

\path
node (up) at (0, -12) {}
node (left) at (-41,-83) {}
node (right) at (41,-83) {};
\draw [->] (up) [one] -- (left);
\draw [->] (left) [one] -- (right);
\draw [->] (up) [one] -- (right);

\path
node (up) at (0, -95) {}
node (left) at (-41,-166) {}
node (right) at (41,-166) {};
\draw [->] (up) [two] -- (left);
\draw [->] (left) [two] -- (right);
\draw [->] (up) [two] -- (right);

\path
node (up) at (113, 71) {}
node (left) at (72,0) {}
node (right) at (154,0) {};
\draw [->] (up) [two] -- (left);
\draw [->] (left) [one] -- (right);
\draw [->] (up) [two] -- (right);

\path
node (up) at (113, -12) {}
node (left) at (72,-83) {}
node (right) at (154,-83) {};
\draw [->] (up) [one] -- (left);
\draw [->] (left) [two] -- (right);
\draw [->] (up) [two] -- (right);

\path
node (up) at (113, -95) {}
node (left) at (72,-166) {}
node (right) at (154,-166) {};
\draw [->] (up) [one] -- (left);
\draw [->] (left) [one] -- (right);
\draw [->] (up) [two] -- (right);

\path
node (up) at (226, 71) {}
node (left) at (185,0) {}
node (right) at (267,0) {};
\draw [->] (up) [two] -- (left);
\draw [->] (left) [zero] -- (right);
\draw [->] (up) [two] -- (right);

\path
node (up) at (226, -12) {}
node (left) at (185,-83) {}
node (right) at (267,-83) {};
\draw [->] (up) [zero] -- (left);
\draw [->] (left) [two] -- (right);
\draw [->] (up) [two] -- (right);

\path
node (up) at (226, -95) {}
node (left) at (185,-166) {}
node (right) at (267,-166) {};
\draw [->] (up) [zero] -- (left);
\draw [->] (left) [two] -- (right);
\draw [->] (up) [zero] -- (right);

\path
node (up) at (339, 71) {}
node (left) at (298,0) {}
node (right) at (380,0) {};
\draw [->] (up) [one] -- (left);
\draw [->] (left) [zero] -- (right);
\draw [->] (up) [one] -- (right);

\path
node (up) at (339, -12) {}
node (left) at (298,-83) {}
node (right) at (380,-83) {};
\draw [->] (up) [one] -- (left);
\draw [->] (left) [zero] -- (right);
\draw [->] (up) [zero] -- (right);

\path
node (up) at (452, 71) {}
node (left) at (411,0) {}
node (right) at (493,0) {};
\draw [->] (up) [zero] -- (left);
\draw [->] (left) [one] -- (right);
\draw [->] (up) [zero] -- (right);

\path
node (up) at (452, -12) {}
node (left) at (411,-83) {}
node (right) at (493,-83) {};
\draw [->] (up) [zero] -- (left);
\draw [->] (left) [one] -- (right);
\draw [->] (up) [one] -- (right);

\path
node (up) at (339, -95) {}
node (left) at (298,-166) {}
node (right) at (380,-166) {};
\draw [->] (up) [one] -- (left);
\draw [->] (left) [zero] -- (right);
\draw [->] (up) [two] -- (right);

\path
node (up) at (452, -95) {}
node (left) at (411,-166) {}
node (right) at (493,-166) {};
\draw [->] (up) [zero] -- (left);
\draw [->] (left) [two] -- (right);
\draw [->] (up) [one] -- (right);

\end{tikzpicture}

\end{definition}

If a triple is not agreeable we call it \emph{disagreeable}. Note that all agreeable triples are \emph{transitive} rather than \emph{cyclic}, i.e. there is both a vertex emitting two arrows and a vertex absorbing two arrows.

\begin{fact}
\label{listfact}
A triple is disagreeable if and only if it is either\dots
\begin{aufz}
\item \dots a cyclic triple, regardless of the colouring, i.e.
\begin{tikzpicture}[line/.style={ultra thick, densely dotted, black}]

\pgftransformscale{0.0037}
\path
node (up) at (0, 362) {}
node (left) at (-209,0) {}
node (right) at (209,0) {};
\draw [->] (up) [line] -- (left);
\draw [->] (left) [line] -- (right);
\draw [->] (right) [line] -- (up);
\end{tikzpicture}

\item \dots or one of the following:

\begin{tikzpicture}[
zero/.style={ultra thick, solid, ZERO},
one/.style={ultra thick, solid, ONE},
two/.style={ultra thick, solid, TWO}
]

\pgftransformscale{0.02}

\path
node (up) at (-113, 71) {}
node (left) at (-154,0) {}
node (right) at (-72,0) {};
\draw [->] (up) [two] -- (left);
\draw [->] (left) [two] -- (right);
\draw [->] (up) [one] -- (right);

\path
node (up) at (0, 71) {}
node (left) at (-41,0) {}
node (right) at (41,0) {};
\draw [->] (up) [two] -- (left);
\draw [->] (left) [one] -- (right);
\draw [->] (up) [one] -- (right);

\path
node (up) at (113, 71) {}
node (left) at (72,0) {}
node (right) at (154,0) {};
\draw [->] (up) [one] -- (left);
\draw [->] (left) [two] -- (right);
\draw [->] (up) [one] -- (right);

\path
node (up) at (-113, -12) {}
node (left) at (-154,-83) {}
node (right) at (-72,-83) {};
\draw [->] (up) [two] -- (left);
\draw [->] (left) [two] -- (right);
\draw [->] (up) [zero] -- (right);

\path
node (up) at (0, -12) {}
node (left) at (-41,-83) {}
node (right) at (41,-83) {};
\draw [->] (up) [two] -- (left);
\draw [->] (left) [zero] -- (right);
\draw [->] (up) [zero] -- (right);

\path
node (up) at (113, -12) {}
node (left) at (72,-83) {}
node (right) at (154,-83) {};
\draw [->] (up) [zero] -- (left);
\draw [->] (left) [zero] -- (right);
\draw [->] (up) [two] -- (right);

\path
node (up) at (226, 71) {}
node (left) at (185,0) {}
node (right) at (267,0) {};
\draw [->] (up) [one] -- (left);
\draw [->] (left) [one] -- (right);
\draw [->] (up) [zero] -- (right);

\path
node (up) at (226, -12) {}
node (left) at (185,-83) {}
node (right) at (267,-83) {};
\draw [->] (up) [zero] -- (left);
\draw [->] (left) [zero] -- (right);
\draw [->] (up) [one] -- (right);

\path
node (up) at (339, 71) {}
node (left) at (298,0) {}
node (right) at (380,0) {};
\draw [->] (up) [two] -- (left);
\draw [->] (left) [one] -- (right);
\draw [->] (up) [zero] -- (right);

\path
node (up) at (339, -12) {}
node (left) at (298,-83) {}
node (right) at (380,-83) {};
\draw [->] (up) [two] -- (left);
\draw [->] (left) [zero] -- (right);
\draw [->] (up) [one] -- (right);

\path
node (up) at (452, 71) {}
node (left) at (411,0) {}
node (right) at (493,0) {};
\draw [->] (up) [one] -- (left);
\draw [->] (left) [two] -- (right);
\draw [->] (up) [zero] -- (right);

\path
node (up) at (452, -12) {}
node (left) at (411,-83) {}
node (right) at (493,-83) {};
\draw [->] (up) [zero] -- (left);
\draw [->] (left) [one] -- (right);
\draw [->] (up) [two] -- (right);

\end{tikzpicture}
\end{aufz}
\end{fact}

We remind the reader of the definition of weakly compacts.

\begin{definition}
$\kappa \in \On \setminus 3$ is called weakly compact iff $\kappa \arrows (\kappa)^2_\omega$.
\end{definition}

\begin{fact}
For all $\kappa \in \On$ we have $\kappa \arrows (\kappa)^2_2$ iff $\kappa \in \{0, 1, 2, \omega\}$ or $\kappa$ is weakly compact.
\end{fact}

\begin{fact}
\label{partitionequivalence}
$\kappa \arrows (\kappa)^2_2$ iff $\forall n < \omega, \alpha < \kappa : \kappa \arrows (\kappa)^n_\alpha$
\end{fact}

Both facts are corollaries of theorem 7.6 of \cite{003K0} which can be found on page 76.

\begin{theorem}
Let $\kappa \in \On \setminus 3$ satisfy $\kappa \arrows (\kappa)^2_2$.

The partition relation $\kappa^2 l \arrows (\kappa^2 m, n)$ holds true if and only if every coloured digraph $C = \langle l, A, c\rangle$ with $\ran(c) = 3$ contains an independent set of size $m$ or there is a subtournament $S$ of $C$ induced by a set of $n$ vertices such that all triples in $S$ are agreeable.
\label{main}
\end{theorem}

\begin{proof}

\begin{aufz}
\item Let us first assume that the finite combinatorial characterization above holds true. Given any colouring $\chi : [\kappa^2 l]^2 \longrightarrow 2$ we define a new colouring $\chi'$ as follows:
\begin{align*}
\chi' : [\kappa]^4 \longrightarrow 2^{4l^2} \\
\{h, j, i, k\}_< \longmapsto \sum_{f, g < l} \big(\? & \chi(\{\kappa^2 f + \kappa h + j, \kappa^2 g + \kappa i + k\}) 2^{4(lf + g)} \\
+ & \chi(\{\kappa^2 f + \kappa h + i, \kappa^2 g + \kappa j + k\}) 2^{4(lf + g) + 1} \\
+ & \chi(\{\kappa^2 f + \kappa h + k, \kappa^2 g + \kappa j + i\}) 2^{4(lf + g) + 2}\\
+ & \chi(\{\kappa^2 f + \kappa h + j, \kappa^2 g + \kappa h + i\}) 2^{4(lf + g) + 3}\?\big)
\end{align*}
Now we use fact \ref{partitionequivalence} in the form $\kappa \arrows (\kappa)^4_{2^{4l^2}}$ thereby finding an $X \in [\kappa]^\kappa$ homogeneous for $\chi'$.

Now fix a surjection $b : X \surj \soast{\langle \nu, \gamma\rangle}{\nu < l \wedge \gamma < \kappa + 1}$ and monotonic enumeration functions $e_\nu : \kappa \longleftrightarrow \pwi{b^{-1}}{(\{\langle\nu, \kappa\rangle\})}$ such that
\begin{aufz}
\item $\forall \nu < l \forall \gamma < \kappa + 1 : \card{\pwi{b^{-1}}{(\{\langle\nu, \gamma\rangle\})}} = \kappa$.
\item $\forall \nu < l \forall \gamma < \kappa \forall \delta \in \pwi{b^{-1}}{(\{\langle\nu, \gamma\rangle\})}: e_\nu(\gamma) < \delta$.
\end{aufz}
Now let $Y := \soast{\kappa^2 f + \kappa h + j}{f < l \wedge h \in b^{-1}(\langle f, \kappa\rangle) \wedge j \in b^{-1}(\langle f, e^{-1}_f(h)\rangle)}$.

In prose what's happening is this. We distribute the elements of $X$ to signify multiples of $1$ or $\kappa$ in the $l$ blocks of size $\kappa^2$. No two ordinals may appear in two roles and within one $\kappa^2$-block the $1$-coordinate is always larger than the $\kappa$-coordinate. Note that $Y \in [\kappa^2 l]^{\kappa^2 l}$.

Now observe that the colour of an element $\{\kappa^2 f + \kappa h + j, \kappa^2 g + \kappa i + k\} \in [Y]^2$---where we suppose that $h < i$---is completely determined by $f$, $g$ and whether $h < j < i < k$, $h < i < j < k$, $h < i < k < j$ or $h = i < j < k$.

Now we define a coloured digraph $C = \langle l, A, c\rangle$ as follows. If for $\{f, g\} \in [l]^2$ we have
\begin{aufz}
\item $\forall \{\kappa^2 f + \kappa h + j, \kappa^2 g + \kappa i + k\} \in [Y]^2 : \chi(\{\kappa^2 f + \kappa h + j, \kappa^2 g + \kappa i + k\}) = 0$ then there is no arc between $f$ and $g$ in $C$.
\item $\exists \{\kappa^2 f + \kappa h + j, \kappa^2 g + \kappa i + k\} \in [Y]^2 : \chi(\{\kappa^2 f + \kappa h + j, \kappa^2 g + \kappa i + k\}) = 1$ where we again assume that $h < i$ and
\begin{aufz}
\item $h < i < k < j$ then we let $f \zbeat g$.
\item $h < i < j < k$ then we let $f \obeat g$.
\item $h < j < i < k$ then we let $f \tbeat g$.
\end{aufz}
In this case there may very well be more than one way to put an arc between $f$ and $g$ and to colour it---clearly there could be six possible ways to do this yet it suffices just to do it in one way.
\end{aufz}
Our remark above shows us that this is well-defined. Since $l$ is natural we do not even have to use the axiom of choice.

Now we may use our combinatorial statement. If $I \in [l]^m$ is independent in $C$ then obviously $Z := \soast{\{\kappa^2 f + \kappa h + j, \kappa^2 g + \kappa i + k\} \in [Y]^2}{\{f, g\} \subset I \wedge h, i, j, k < \kappa}$ is a homogeneous set of size $\kappa^2 m$ in colour $0$. So let us assume that there is a subtournament $S$ of $C$ of size $n$ such that all triples in $S$ are agreeable.

Let $\langle m_0, m_2, \dots, m_{n-1}\rangle$ be the sequence of natural numbers defined by the following backwards recursion:
\begin{aufz}
\item $m_n := 0$,
\item $m_{i-1} := 3m_i + 2$.
\end{aufz}
We have $m_i:=3^{n-i}-1$. Furthermore let $\langle v_0, v_1, \dots, v_{n-1}\rangle$ be a sequence of vertices such that $S := \soast{v_i}{i < n}$.

Now we are going to inductively construct the $1$-homogeneous set $\soast{\kappa a_i + b_i}{i < n}$ where $a_i$ and $b_i$ are ordinals less than $\kappa$ for $i < n$.

First we are going to inductively define a sequence $\langle s_i | i < m_0\rangle$ of ordinals less than $\kappa$. Let $s_0 := 0$. Given $\langle s_j |j \leqslant i\rangle$ choose $\alpha_i < \kappa$ such that $e_{v_j}(\alpha_i) \geqslant s_i$ for all $j < n$ and then choose $s_{i+1} < \omega$ such that $\forall j < n, k \leqslant i : \pwi{b^{-1}}{(\{\langle v_j, \alpha_k \rangle\})} \cap s_{i+1} \setminus s_i \supsetneq \emptyset$.

In step $i < n$ suppose we have already defined a $1$-homogeneous set of size $i$.

The following will be our induction hypothesis.
\begin{enumerate}
\item $\min\soast{j \in \omega \setminus 1}{\exists k < \omega : \soast{a_h, b_h}{h < i} \cap s_{j+k+1} \setminus s_{j+k} \supsetneq \emptyset \wedge \soast{a_h, b_h}{h < i} \cap s_{k+1} \setminus s_k \supsetneq \emptyset}>m_{i+1}$,
\item $\forall j, k < i : v_j \zbeat v_k \leftrightarrow a_j < a_k < b_k < b_j$,
\item $\forall j, k < i : v_j \obeat v_k \leftrightarrow a_j < a_k < b_j < b_k$,
\item $\forall j, k < i : v_j \tbeat v_k \leftrightarrow a_j < b_j < a_k < b_k$.
\end{enumerate}
This we do as follows, let $L_\kappa := \max(\soast{a_h}{h < i \wedge v_h \beat v_i} \cup \soast{b_h}{h < i \wedge v_h \tbeat v_i})$ and $U_\kappa := \min(\soast{a_h}{h < i \wedge v_i \beat v_h} \cup \soast{b_h}{h < i \wedge v_h \zbeat v_i} \cup \soast{b_h}{h < i \wedge v_h \obeat v_i})$.

\begin{claim}
\label{firstclaim}
$L_\kappa < U_\kappa$.
\end{claim}
\begin{nestedproof}

Suppose not. The claim could fail in the following ways:

\begin{aufz}
\item $\max\soast{a_h}{h < i \wedge v_h \beat v_i}\geqslant \min\soast{a_j}{j < i \wedge v_i\beat v_j}$. Fix $h$ and $j$ witnessing this. The induction hypothesis implies that $v_j \beat v_h$. But then $v_h$, $v_i$ and $v_j$ form a cyclic and hence disagreeable triple. \con
\item $\max\soast{a_h}{h < i \wedge v_h \beat v_i}\geqslant\min\soast{b_j}{j < i \wedge v_j \zbeat v_i}$. As before fix witnesses $h$ and $j$. By induction hypothesis we have $v_j \tbeat v_h$ and fact \ref{listfact} implies that $\{v_h, v_i, v_j\}$ is disagreeable. \con
\item $\max\soast{a_h}{h < i \wedge v_h \beat v_i}\geqslant\min\soast{b_j}{j < i \wedge v_j \obeat v_i}$. Fix $h$ and $j$ witnessing this. The induction hypothesis implies that $v_j \tbeat v_h$. A look at fact \ref{listfact} shows that $\{v_h, v_i, v_j\}$ is disagreeable. \con
\item $\max\soast{b_h}{h < i \wedge v_h \zbeat v_i}\geqslant\min\soast{a_j}{j < i \wedge v_i\beat v_j}$. Fix $h, j$ witnessing this. By induction hypothesis we now either have $v_j\beat v_h$, $v_h \zbeat v_j$ or $v_h \obeat v_j$. In every case $\{v_h, v_i, v_j\}$ is disagreeable. \con
\item $\max\soast{b_h}{h < i \wedge v_h \tbeat v_i}\geqslant\min\soast{b_j}{h < i \wedge v_j \zbeat v_i}$. Let $h, j$ witness this. As before the induction hypothesis implies that either $v_h \zbeat v_j$, $v_j \obeat v_h$ or $v_j \tbeat v_h$. Again by fact \ref{listfact} $\{v_h, v_i, v_j\}$ is disagreeable. \con
\item $\max\soast{b_h}{h < i \wedge v_h \tbeat v_i}\geqslant\min\soast{b_j}{h < i \wedge v_j \obeat v_i}$. Fix witnesses $h$ and $j$ to this fact. By induction hypothesis we now either have $v_h \zbeat v_j$, $v_j \obeat v_h$ or $v_j \tbeat v_h$. Fact \ref{listfact} again shows $\{v_h, v_i, v_j\}$ to be disagreeable. \con
\end{aufz}
\end{nestedproof}

Let $h < m_0$ be minimal such that $L_\kappa < s_h$. Then take

\begin{align*}
a_i \in \pwi{e_{v_i}}{(\kappa)} \cap s_{h + m_i + 1} \setminus s_{h + m_i}.
\end{align*}

Now we are going to define $b_i$. Similarly, let $L_1 := \max(\{a_i\} \cup \soast{a_h}{h < i \wedge v_i \obeat v_h} \cup \soast{b_h}{h < i \wedge v_i \zbeat v_h} \cup \soast{b_h}{h < i \wedge v_h \obeat v_i})$ and $U_1 := \min(\soast{a_h}{h < i \wedge v_i\tbeat v_h} \cup \soast{b_h}{h < i \wedge v_h\zbeat v_i} \cup \soast{b_h}{h < i \wedge v_i\obeat v_h})$.

\begin{claim}
\label{secondclaim}
$L_1 < U_1$.
\end{claim}

\begin{nestedproof}
As above we distinguish cases and reach contradictions.

\begin{aufz}
\item Suppose $a_i \geqslant \min\soast{a_h}{h < i \wedge v_i \tbeat v_h}$. Since $\soast{a_h}{h < i \wedge v_i \tbeat v_h} \subset \soast{a_h}{h < i \wedge v_i \beat v_h}$ we have $\min\soast{a_h}{h < i \wedge v_i \tbeat v_h} \geqslant \min\soast{a_h}{h < i \wedge v_i \beat v_h}$ and since $\min\soast{a_h}{h < i \wedge v_i \beat v_h} \geqslant U_\kappa$ it follows that $a_i \geqslant U_\kappa$. \con
\item If $a_i \geqslant \min\soast{b_h}{h < i \wedge v_h \zbeat v_i}$ then analogously we get $\min\soast{b_i}{h < i \wedge v_h \zbeat v_i} \geqslant U_\kappa$ implying $a_i \geqslant U_\kappa$. \con
\item Suppose $a_i \geqslant \min\soast{b_h}{h < i \wedge v_i \obeat v_h}$. Similarly we have $\min\soast{b_h}{h < i \wedge v_h \obeat v_i} \geqslant U_\kappa$ and thus $a_i \geqslant U_\kappa$. \con
\item Suppose $\max\soast{a_h}{h < i \wedge v_i \obeat v_h} \geqslant \min\soast{a_j}{j < i \wedge v_i \tbeat v_j}$. Fix witnesses $h$ and $j$ to this fact. The induction hypothesis implies $v_j \beat v_h$. Fact \ref{listfact} shows that in all three cases $\{v_h, v_i, v_j\}$ is disagreeable.
\item Suppose $\max\soast{a_h}{h < i \wedge v_i \obeat v_h} \geqslant \min\soast{b_j}{j < i \wedge v_j \zbeat v_i}$. Once more fix witnesses $h, j$. Again the induction hypothesis implies that $v_j \tbeat v_h$. Fact \ref{listfact} implies that $\{v_h, v_i, v_j\}$ is a disagreeable triple. \con
\item In case $\max\soast{a_h}{h < i \wedge v_i \obeat v_h} \geqslant \min\soast{b_j}{j < i \wedge v_i\obeat v_j}$ we can again fix $h$ and $j$ witnessing this and use the induction hypothesis to see that $v_j \tbeat v_h$. A backcheck with fact \ref{listfact} assures us of the disagreeability of $\{v_h, v_i, v_j\}$. \con
\item In case $\max\soast{b_h}{h < i \wedge v_i \zbeat v_h} \geqslant \min\soast{a_j}{j < i \wedge v_i \tbeat v_j}$ let us once more fix witnesses $h$ and $j$ to this fact and use the induction hypothesis. We must either have $v_j \beat v_h$ or $v_h \zbeat v_j$ or $v_h \obeat v_j$. Fact \ref{listfact} shows that in all five cases $\{v_h, v_i, v_j\}$ is disagreeable. \con
\item Suppose that $\max\soast{b_h}{h < i \wedge v_i \zbeat v_h} \geqslant \min\soast{b_j}{j < i \wedge v_j \zbeat v_i}$. Fix witnesses $h$ and $j$ and use the induction hypothesis in order to see that we either have $v_h \zbeat v_j$ and hence a cyclic triple $\{v_h, v_i, v_j\}$ or we have either $v_j \obeat v_h$ or $v_j \tbeat v_h$. In both of the last two cases $\{v_h, v_i, v_j\}$ is disagreeable by fact \ref{listfact}. \con
\item If $\max\soast{b_h}{h < i \wedge v_i \zbeat v_h} \geqslant\min\soast{b_j}{j < i \wedge v_i \obeat v_j}$ we again fix witnesses $h$ and $j$. By induction hypothesis we have $v_h \zbeat v_j$ or $v_j \obeat v_h$ or $v_j \tbeat v_h$. Fact \ref{listfact} assures us once again of $\{v_h, v_i, v_j\}$'s disagreeability in any case. \con
\item Suppose $\max\soast{b_h}{h < i \wedge v_h \obeat v_i} \geqslant \min\soast{b_j}{j < i \wedge v_j \zbeat v_i}$. Let $h$ and $j$ be witnesses to this fact. The induction hypothesis implies that either $v_h \zbeat v_j$ or $v_j \obeat v_h$ or $v_j \tbeat v_h$. By fact \ref{listfact} in any case $\{v_h, v_i, v_j\}$ is disagreeable. \con
\item If $\max\soast{b_h}{h < i \wedge v_h \obeat v_i} \geqslant \min\soast{b_j}{j < i \wedge v_i \obeat v_j}$ we fix witnesses $h$ and $j$. By induction hypothesis we either have $v_h \zbeat v_j$ or $v_j \obeat v_h$ or $v_j \tbeat v_h$. In each case by fact \ref{listfact} $\{v_h, v_i, v_j\}$ would be disagreeable. \con
\item If $\max\soast{b_h}{h < i \wedge v_h \obeat v_i} \geqslant \min\soast{a_j}{j < i \wedge v_i \tbeat v_j}$ we can again fix witnesses $h, j$ to this fact. By induction hypothesis we can conclude that either $v_j \beat v_h$ in which case $\{v_h, v_i, v_j\}$ is cyclic and in particular disagreeable or that either $v_h \zbeat v_j$ or $v_h \obeat v_j$. But fact \ref{listfact} tells us that in these cases $\{v_h, v_i, v_j\}$ is disagreeable too. \con
\end{aufz}
\end{nestedproof}

Now let $h < m_0$ be minimal such that $L_1 < s_h$. Take
\begin{align*}
b_i \in \pwi{b^{-1}}{(\{\langle v_i, e^{-1}(a_i)\rangle\})} \cap s_{h + 2m_i + 1} \setminus s_{h + 2m_i}\footnotemark.
\end{align*}

By claims \ref{firstclaim} and \ref{secondclaim} and condition 1. of the inductive hypothesis we have chosen $a_i$ and $b_i$ such that $a_i < U_\kappa$ and $b_i < U_1$. So by definition of the bounds $L_\kappa$, $U_\kappa$, $L_1$\ox\ and $U_1$ conditions 2. to 4. of our inductive hypothesis are fulfilled. We also made sure that 1. holds true for the next step.

But then, after $n$ steps, by the definition of our coloured digraph $C$, we see that $\soast{\kappa a_i + b_i}{i < n}$ is homogeneous of colour $1$.

\footnotetext{There are basically two ideas behind this abundance of case distinctions. First one can conceive of a completeness theorem for linear orders. Unless one cannot derive a contradiction there has to be a model. The other is the fact that if a tournament has a cycle it has one of length $3$. Here it only gets a little bit more involved since one has to simultaneously consider the $\omega$-coordinates and the $1$-coordinates.}

\item Now suppose towards a contradiction that there is a a $3$-coloured digraph $D = \langle l, A, c\rangle$ of size $l$ such that both all its independent sets have size less than $m$ and every subtournament $S$ of $D$ of size $n$ contains a disagreeable triple.

Now we define a colouring $\chi$ as follows.

\begin{align*}
\chi : [\omega^2 l] ^2&\longrightarrow 2\\
\{\omega^2 e + \omega g + i, \omega^2 f + \omega h + j\} &\longmapsto \begin{cases}
1 \text{ if } g < h < j < i \text{ and } e \zbeat f,\\
1 \text{ if } g < h < i < j \text{ and } e \obeat f,\\
1 \text{ if } g < i < h < j \text{ and } e \tbeat f,\\
0 \text{ if } e \text{ and } f \text{ are unconnected,}\\
0 \text{ if } \card{\{g, h, i, j\}} < 4,\\
0 \text{ if } e = f \text{ or } i < g \text{ or } j < h.
\end{cases}
\end{align*}

In the naming of the pair of ordinals we can of course assume without loss of generality that $g \leqslant h$.

Now if there were a set $H \in [\kappa^2 l]^{\kappa^2 m}$ which is $0$-homogeneous for $\chi$ there would be a set $S \in [l]^m$ such that for all $s \in S$ we have $\otyp(\soast{\kappa^2 s + \alpha \in H}{\alpha < \kappa^2}) = \kappa^2$. Since all independent sets in $D$ have size less than $m$ it follows that we have $e \beat f$ for $e, f \in S$. Now if $e \tbeat f$, choose first $g \in \soast{k}{k < \kappa \wedge \card{\soast{k'}{k' < \kappa \wedge \kappa^2 e + \kappa k + k' \in H}} = \kappa}$, then $i \in \soast{k}{k \in \kappa \setminus (g + 1) \wedge \kappa^2 e + \kappa g + k \in H}$, then $h \in \soast{k}{k \in \kappa \setminus (i + 1) \wedge \card{\soast{k'}{k' < \kappa \wedge \kappa^2 f + \kappa k + k' \in H}} = \kappa}$ and finally $j \in \soast{k}{k \in \kappa \setminus (h + 1) \wedge \kappa^2 f + \kappa h + k \in H}$. But now we have found $\kappa^2 e + \kappa g + i, \kappa^2 f + \kappa h + j \in H$ such that $\chi(\{\kappa^2 e + \kappa g + i, \kappa^2 f + \kappa h + j\}) = 1$ which contradicts our assumption that $H$ was $0$-homogeneous. If instead of $e \tbeat f$ we have $e \obeat f$ or $e \zbeat f$ we can argue similarly, we only have to change the sequence of choices.

So assume that there is a set $H \in [\omega^2 l]^n$ which is $1$-homogeneous for $\chi$. Considering the last clause in the definition of $\chi$ we conclude that no $\kappa^2$-block can contain more than one element of $H$, i.e. we have $\card{H \cap \kappa^2(k + 1) \setminus \kappa^2 k} < 2$ for every $k < l$. So $S := \soast{s}{s < l \wedge H \cap \kappa^2(s + 1) \setminus \kappa^2 s \supsetneq \emptyset}$ has size $n$. Moreover the fourth clause in the definition of $\chi$ implies that every two elements of $S$ are connected by an arrow. So $S$ spans a subtournament of $D$. Hence $S$ has to contain a disagreeable triple $T = \{t_0, t_1, t_2\} \in [S]^3$. Take $\kappa^2 t_0 + \kappa u_0 + v_0, \kappa^2 t_1 + \kappa u_1 + v_1, \kappa^2 t_2 + \kappa u_2 + v_2 \in H$. Now we distinguish four cases

\begin{aufz}

\item $T$ is cyclic. Suppose we have $t_0 \beat t_1 \beat t_2 \beat t_0$. Considering the first three clauses of the definition of $\chi$ we may conclude that for the elements chosen above we have $u_0 < u_1 < u_2 < u_0$. \con

\item Suppose $t_0 \obeat t_1$ or $t_0 \tbeat t_1$, that $t_2 \zbeat t_1$ or $t_1 \obeat t_2$ or $t_1 \tbeat t_2$\ox\ and that $t_0 \zbeat t_2$. Then using again the definition of $\chi$ we arrive at $v_0 < v_1 < v_2 < v_0$. \con

\item Suppose that $t_0 \tbeat t_1$, that $t_1 \beat t_2$.\ox\ and that $t_0 \zbeat t_2$ or $t_0 \obeat t_2$. Then $v_0 < u_1 < u_2 < v_0$ follows. \con

\item Suppose that $t_0 \tbeat t_1$, that $t_2 \zbeat t_1$ or $t_2 \obeat t_1$\ox\ and that either $t_0 \zbeat t_2$ or $t_2 \obeat t_0$ or $t_2 \tbeat t_0$. Then $v_0 < u_1 < v_2 < v_0$. \con

\end{aufz}

Note that modulo a renaming of $t_0$, $t_1$\ox\ and $t_2$ the above clauses cover all disagreeable triples. So this finishes the second part of the proof.

\end{aufz}

\end{proof}

The author wants to point out that if you believe in enough reflection you might believe something like (There are stationarily many weakly compact cardinals in the universe.) so theorem \ref{main} really says something about a lot of ordinals.

So we now can also state---not in complete earnest---a joint generalization of the theorems \ref{erdosradomodern} and \ref{main}. In this artificial sense being agreeable is a generalization of being transitive.

\begin{theorem}
For $i \in 3 \setminus 1$ and $\kappa \in \On \setminus 3$ such that $\kappa \arrows (\kappa)^2_2$ the partition relation $\kappa^i l \arrows (\kappa^i m, n)$ holds true if and only if every coloured digraph $C = \langle l, A, c\rangle$ with $\ran(c) = 2i - 1$ contains an independent set of size $m$ or there is a
subtournament $S$ of $C$ induced by a set of $n$ vertices such that all triples in $S$ are agreeable.
\end{theorem}

The special case the author was actually most interested in however is this one:

\begin{corollary}

The partition relation $\omega^2 l \arrows (\omega^2 m, n)$ holds true if and only if every coloured digraph $C = \langle l, A, c\rangle$ with $\ran(c) = 3$ contains an independent set of size $m$ or there is a subtournament $S$ of $C$ induced by a set of $n$ vertices such that all triples in $S$ are agreeable.

\end{corollary}

Henceforth $A_n$ is a graph on $n$ vertices such that all its induced subgraphs on three vertices are agreeable. This allows us to define the Ramsey number $r(I_m, A_n)$ to be the smallest natural number such that any digraph on $r(I_m, A_n)$ vertices either contains an independent set of size $m$ or an induced subgraph on $n$ vertices all triples of which are agreeable. Since $A_n$ does not denote a single graph but a class of graphs this might be seen as a slight abuse of notation, or at least a shorthand. So the equivalence between the infinitary partition problems and finite digraph combinatorics above can be also condensed into the following formula:

\begin{align*}
r(\kappa^2 m, n) = \kappa^2 r(I_m, A_n) \text{ for }\kappa = \omega \text{ or } \kappa \text{ weakly compact.}
\end{align*}

\begin{theorem}
$r(I_2, A_3) > 9$.
\end{theorem}

\begin{proof}
The statement above holds because there is the following tricoloured nine-person tournament:

\begin{tikzpicture}[scale=2,
two/.style={ultra thick, solid, TWO},
one/.style={ultra thick, solid, ONE},
zero/.style={ultra thick, solid, ZERO},
vertex/.style={ultra thick, circle, draw=dark grey, fill=very light grey, minimum size=3cc}]

\pgftransformscale{0.95}

\draw (0,0) arc (0:40:3) [white] node [vertex] (a) {}

arc (40:80:3) node (b) [vertex] {}

arc (80:120:3) node (c) [vertex] {}

arc (120:160:3) node (d) [vertex] {}

arc (160:200:3) node (e) [vertex] {}

arc (200:240:3) node (f) [vertex] {}

arc (240:280:3) node (g) [vertex] {}

arc (280:320:3) node (h) [vertex] {}

arc (320:360:3) node (i) [vertex] {};

\foreach \from/\to in {a/b,b/c,c/d,d/e,e/f,f/g,g/h,h/i,i/a}
\draw [->] (\from) [two] -- (\to);

\foreach \from/\to in {a/c,b/d,c/e,d/f,e/g,f/h,g/i,h/a,i/b, a/d,b/e,c/f,d/g,e/h,f/i,g/a,h/b,i/c}
\draw [->] (\from) [one] -- (\to);

\foreach \from/\to in {a/e,b/f,c/g,d/h,e/i,f/a,g/b,h/c,i/d}
\draw [->] (\from) [zero] -- (\to);

\end{tikzpicture}

Algebraically this can be seen as the tournament on $\mathbb{Z}_9$ where from every $i \in \mathbb{Z}_9$ there is a \two\ arrow pointing at $i + 1$, \one\ arrows pointing both at $i + 2$ and $i + 3$\ox\ and a \zero\ arrow pointing at $i + 4$.

A close observation reveals that every triple in this tournament is disagreeable thereby showing $\omega^2 9 \doesntarrow (\omega^2 2, 3)$.

There is even an easier counterexample which is this one:

\begin{tikzpicture}[scale=2,
two/.style={ultra thick, solid, TWO},
one/.style={ultra thick, solid, ONE},
vertex/.style={ultra thick, circle, draw=dark grey, fill=very light grey, minimum size=3cc}]

\pgftransformscale{0.95}

\draw (0,0) arc (0:40:3) [white] node [vertex] (a) {}

arc (40:80:3) node (b) [vertex] {}

arc (80:120:3) node (c) [vertex] {}

arc (120:160:3) node (d) [vertex] {}

arc (160:200:3) node (e) [vertex] {}

arc (200:240:3) node (f) [vertex] {}

arc (240:280:3) node (g) [vertex] {}

arc (280:320:3) node (h) [vertex] {}

arc (320:360:3) node (i) [vertex] {};

\foreach \from/\to in {a/d,a/e,a/f,b/d,b/e,b/f,c/d,c/e,c/f, d/g,d/h,d/i,e/g,e/h,e/i,f/g,f/h,f/i, g/a,g/b,g/c,h/a,h/b,h/c,i/a,i/b,i/c}
\draw [->] (\from) [one] -- (\to);

\foreach \from/\to in {a/b,b/c,c/a,d/e,e/f,f/d,g/h,h/i,i/g}
\draw [->] (\from) [two] -- (\to);

\end{tikzpicture}

This is algebraically interpretable as the tournament on $\mathbb{Z}_3 \times \mathbb{Z}_3$ where from every $\langle i, j\rangle$ there is a \two\ arrow pointing at $\langle i + 1, j\rangle$ and \one\ arrows pointing at $\langle i, j + 1\rangle$ and $\langle i + 1, j + 1\rangle$ and $\langle i + 2, j + 1\rangle$.

Also this should be easily checkable. Either a triple is one of three completely coloured in \two---these are disagreeable---or it shares exactly one point with each of these triples---then it forms a disagreeable triple completely coloured in \one---or it shares two points with one and one with another---then it is either one of two \two/\one-coloured disagreeable triples
. These counterexamples were not found in a very systematic way, so there might be others or there might not be others.
\end{proof}

Now some easy observations quickly yield upper bounds for $r(I_2, A_3)$. For example by ordering the vertices and then colour a pair of them in one of six colours depending on both whether the direction of the arc agrees with the ordering or not and the original colour of the arc. Since all monochromatic transitive triples are agreeable this yields $r(I_2, A_3) \leqslant r(3, 3, 3, 3, 3, 3)$. Using $r(3, 3, 3) = 17$ from \cite{955GG0} together with the recursive formula we get $r(I_2, A_3) \leqslant 1898$.

A much better bound is found by observing that every four-person-\linebreak
tournament contains a transitive triple. So clearly $r(I_2, A_3) \leqslant r(4, 4, 4)$. The recursive formula together with the result that $r(3, 3, 4) \leqslant 31$ from \cite{998PR0} yields $r(4, 4, 4) \leqslant 236$ so $r(I_2, A_3) \leqslant 236$ follows. Finally by observing that every six-person-tournament coloured in either \zero\ or \one\ contains an agreeable triple we get $r(I_2, A_3) \leqslant r(4, 6)$. In \cite{997MR0} it is shown that $r(4, 6) \leqslant 41$ so we have $r(I_2, A_3) \leqslant 41$.

But far better results are possible.

\begin{fact}
In a tournament on ten vertices containing no agreeable triple there can be at most ten \two\ arrows, thirty \one\ arrows\ox\ and ten \zero\ arrows. So there are at least five \two , twenty-five \one\ox\ and five \zero\ arrows.
\label{arcnum}
\end{fact}

\begin{proof}
Take any tournament on ten vertices not containing any agreeable triple.

Suppose that there are eleven \two\ arrows. By the pigeonhole principle there has to be a vertex where two \two\ arrows leave. But then this vertex together with the targets of the two \two\ arrows forms an agreeable triple.

Suppose that there are eleven \zero\ arrows. Again by the pigeonhole principle there has to be a vertex where two \zero\ arrows leave. Analogously this vertex together with the targets of the two \zero\ arrows forms an agreeable triple.

Suppose that there are thirty-one \one\ arrows. By the pigeonhole principle there has to be a vertex where four \one\ arrows leave. Then all of the targets of these arrows have to be connected by \two\ arrows to one another. Fix one of these target-vertices $v$. Since it is adjacent to at least three \two\ arrows at least two \two\ arrows are leaving from there or two \two\ arrows are arriving there. Suppose two \two\ arrows are leaving. Then $v$ together with the targets of these arrows forms an agreeable triple. If on the other hand two \two\ arrows are arriving then again $v$ together with the sources of these arrows forms an agreeable triple.

The second statement follows immediately after one observes that a tournament on ten vertices has $\binom{10}{2}$, i.e. $45$ arrows.
\end{proof}

\begin{lemma}
Any tournament $T$ on ten vertices containing no agreeable triple contains three \two\ triples, i.e. we have $v_0, \dots, v_9 \in V_T$ with $v_1 \tbeat v_2 \tbeat v_3 \tbeat v_1$ and $v_4 \tbeat v_5 \tbeat v_6 \tbeat v_4$ and $v_7 \tbeat v_8 \tbeat v_9 \tbeat v_7$.
\label{threetriples}
\end{lemma}

\begin{proof}
First observe that any such tournament has to contain two such \two\ triples. This is because by fact \ref{arcnum} there have to be at least twenty-five \one\ arrows and hence five vertices where three \one\ arrows leave. Since the targets of the three \one\ arrows leaving such a vertex have to form a \two\ triple and since any vertex in such a \two\ triple may be the target of at most three \one\ arrows the pigeonhole principle implies that there have to be at least two \two\ triples.

Now suppose towards a contradiction that there are exactly two such \two\ triples $\{v_0, v_1, v_2\}$ and $\{v_3, v_4, v_5\}$ and that there are six vertices where three \one\ arrows leave. Again the targets of these \one\ arrows have to form \two\ triples. Since two \two\ triples cannot overlap without producing an agreeable triple we can apply the pigeonhole principle and conclude that each vertex in each of both triples has to receive three \one\ arrows. Now since three \one\ arrows each are hitting $v_0$, $v_1$\ox\ and $v_2$ the sources of them again form a \two\ triple. Since there are only two of them and $v_0$, $v_1$\ox\ and $v_2$ are connected to each other by \two\ arrows we know that the sources of them can only be $v_3$, $v_4$\ox\ and $v_5$. So we have $v_i \obeat v_j$ for every $i \in 6 \setminus 3$ and $j \in 3$. But since $v_3$ is also hit by three \one\ arrows and is already connected in the opposite direction with $v_0$, $v_1$\ox\ and $v_2$ there has to be a third \two\ triple.

An analogous argument works for the case that there are six vertices where three \one\ arrows \emph{leave}.

So there are at most five vertices where three \one\ arrows leave and at most five vertices where three \one\ arrows arrive. So there are at most twenty-five \one\ arrows. By fact \ref{arcnum} we may restrict our attention to the case of ten \two, twenty-five \one\ox\ and ten \zero\ arrows. Since no vertex is the target of more than one \two\ arrow and no vertex is the source of more than one \two\ arrow the remaining four \two\ arrows have to form a cycle. So without loss of generality we have $v_6 \tbeat v_7 \tbeat v_8 \tbeat v_9 \tbeat v_6$.

Now we argue as before. Consider now again the five vertices where three \one\ arrows are leaving. The targets of the three \one\ arrows leaving such a vertex must form a \two\ triple. By the pigeonhole principle every vertex in one of the two \two\ triples has to receive three \one\ arrows each. But then their sources again form a \two\ triple. This shows that without loss of generality we have $v_i \obeat v_j$ for every $i \in 6 \setminus 3$ and $j \in 3$.

In total there are five vertices with three \one\ arrows leaving and there are five vertices with three \one\ arrows arriving. Both the targets of the former and the sources of the latter have to form \two\ triples. Since every vertex may emit or absorb at most three \one\ arrows it follows that there are two-element-sets $\{w_0, w_1\}$ and $\{w_2, w_3\}$ with $w_0 \tbeat w_1$ and $w_2 \tbeat w_3$ such that $v_i \obeat w_j$ for every $i \in 3$ and $j \in 2$ and $w_i \obeat v_j$ for every $i \in 4 \setminus 2$ and $j \in 6 \setminus 3$.

Now we are going to distinguish some---not altogether different---cases. Note that since there are ten \zero\ arrows and a vertex can emit at most one, every vertex has to emit exactly one.

\begin{aufz}
\item $\card{\{w_0, w_1, w_2, w_3\}} = 2$. Since there are only two possible targets---the elements of the set $\{v_6, \dots, v_9\} \setminus \{w_0, w_1\}$---for the \zero\ arrows emitted by the vertices $v_0, \dots, v_5, w_0, w_1$ one of them has to be target of at least four \zero\ arrows which contradicts fact \ref{arcnum}.

\item $\card{\{w_0, w_1, w_2, w_3\}} = 3$. Let $X := \{w_0, w_1\} \cap \{w_2, w_3\}$. Obviously $\card{X} = 1$. Let $w := \bigcup X$. The only possible target for the \zero\ arrow leaving $w$ is $v := \bigcup(\{v_6, \dots, v_9\} \setminus \{w_0, w_1, w_2, w_3\})$ since for each $i \in 4$ we either have $w = w_i$ or $w \tbeat w_i$ or $w_i \tbeat w$ and for all $i \in 6$ we either have $w \obeat v_i$ or $v_i \obeat w$. Now since the vertices $v_0, \dots, v_5$ all also have to emit a \zero\ arrow and the sources of \zero\ arrows pointing at the same vertex have to be connected by \two\ arrows, lest there be an agreeable triple, we know that none of these \zero\ arrows can point at $v$. None can point at $w$ either since $w$ was already shown to be connected differently to every other vertex. So let $x := \bigcup(\{w_0, w_1\} \setminus \{v, w\})$ and $y := \bigcup(\{w_2, w_3\}\setminus \{v, w\})$. Since $v_i \obeat x$ for any $i \in 3$ it follows that $v_i \zbeat y$ for any $i \in 3$ Similarly $y \obeat v_i$ for any $i \in 6\setminus 3$ implies $v_i \zbeat x$ for any $i \in 6\setminus 3$. But now we see that not both $x$ and $y$ can emit a \zero\ arrow hence we have a contradiction.

\item $\card{\{w_0, w_1, w_2, w_3\}} = 4$. Remember that $v_6, \dots, v_9$ form a \two\ four-cycle. Since $\{v_6, \dots, v_9\} = \{w_0, \dots, w_3\}$ and $w_0 \zbeat w_1$ and $w_2 \tbeat w_3$ we can conclude $w_0 \tbeat w_1 \tbeat w_2 \tbeat w_3 \tbeat w_0$. Now recall that we focused on the case where there are five vertices emitting three \one\ arrows. Since in total there are twenty-five \one\ arrows every arrow has to emit at least two \one\ arrows. Now consider the vertex $w_0$. Since it can only send one \one\ arrow to $w_2$ and none to $v_0, v_1, v_2, w_1$ or $w_3$ it has to send one to $v_3, v_4$ or $v_5$. But then $w_0$ and $w_2$ would have to be connected by a \two\ arrow making either $\{w_0, w_1, w_2\}$ or $\{w_2, w_3, w_0\}$ an agreeable triple. Again we have a contradiction.
\end{aufz}
\end{proof}

\begin{theorem}
$r(I_2, A_3) \leqslant 10$.
\label{firstnontrivial}
\end{theorem}

\begin{proof}
Suppose towards a contradiction that there is a tricoloured tournament on ten vertices $v_0, \dots, v_9$ not containing any agreeable triple. By lemma \ref{threetriples} we know that there are three \two\ triples $\{v_1, v_2, v_3\}$, $\{v_4, v_5, v_6\}$ and $\{v_7, v_8, v_9\}$. Since no vertex can emit more than one \two\ arrow and no vertex can absorb more than one \two\ arrow either, we know that $v_0$ has to be connected to $v_1, \dots, v_9$ purely by \zero\ and \one\ arrows. In the situation in which all triples are agreeable we know that the targets of three \one\ arrows starting at $v_0$ have to form a \two\ triple. The same holds for the sources of three \one\ arrows ending at $v_0$ and the sources of three \zero\ arrows ending at $v_0$. Moreover we know that $v_0$ may emit at most one \zero\ arrow. So we may distinguish three cases up to isomorphism:

\begin{aufz}
\item $v_1, v_2, v_3 \zbeat v_0$,\\
$v_4, v_5, v_6 \obeat v_0$\\
$v_0 \obeat v_7, v_8$ and $v_0 \zbeat v_9$.\\
In this case $\{v_0, v_7, v_9\}$ is agreeable.
\item $v_1, v_2, v_3 \zbeat v_0$,\\
$v_0 \obeat v_4, v_5, v_6$\\
$v_7, v_8 \obeat v_0$ and $v_0 \zbeat v_9$.\\
Here $\{v_0, v_8, v_9\}$ is agreeable.
\item $v_1, v_2, v_3 \obeat v_0$,\\
$v_0 \obeat v_4, v_5, v_6$,\\
$v_7, v_8 \zbeat v_0$( and either $v_9 \zbeat v_0$ or $v_0 \zbeat v_9$).
\label{nontrivialcase}

Now we prove two claims:

\begin{claim}
The \two\ triple $\{v_7, v_8, v_9\}$ is connected purely by \one\ arrows to one of the other two \two\ triples.
\label{purelyone}
\end{claim}

\begin{nestedproof}
First note that no \one\ arrow can point from a $v_i$ with $i \in 4\setminus 1$ to a $v_j$ with $j \in 10\setminus 7$, otherwise $\{v_0, v_i, v_j\}$ would be agreeable. Similarly no \one\ arrow may point from a $v_i$ with $i \in 10\setminus 7$ to a $v_j$ with $j \in 7\setminus 4$, otherwise $\{v_0, v_i, v_j\}$ would be agreeable. Since at least two \zero\ arrows are leaving the set $\{v_7, v_8, v_9\}$ heading towards $v_0$ and there can be at most one \zero\ arrow leaving a vertex we conclude that there can be at most one \zero\ arrow leaving $\{v_7, v_8, v_9\}$ in the direction of any vertex $v_i$ with $i \in 10\setminus 1$. So we have either that all arrows between $\{v_7, v_8, v_9\}$ and $\{v_4, v_5, v_6\}$ are either incoming \one\ or incoming \zero\ or that all arrows between $\{v_7, v_8, v_9\}$ and $\{v_1, v_2, v_3\}$ are either outgoing \one\ or incoming \zero. Suppose the latter happens. Assume towards a contradiction that $v_i \zbeat v_j$ for $i \in 4\setminus 1$ and $j \in 10\setminus 7$. Since at most one \zero\ arrow may leave a vertex we know that $v_k \obeat v_i$ for both $k \in 10\setminus(7 \cup \{j\})$. Now choose $k \in 10\setminus(7 \cup\{j\})$ such that $v_k \tbeat v_j$. Then $\{v_i, v_j, v_k\}$ is agreeable, contradiction! So now assume that the former happens, i.e. that all arrows between $\{v_7, v_8, v_9\}$ and $\{v_4, v_5, v_6\}$ are either incoming \one\ or incoming \zero. Assume towards a contradiction that for some $i \in 7\setminus 4$ and some $j \in 10\setminus 7$ we have $v_i \zbeat v_j$. Since at most one \zero\ arrow may leave any vertex it follows that $v_i \obeat v_k$ for $k \in 10\setminus(7\cup\{j\})$. Choose $k \in 10\setminus(7\cup\{j\})$ such that $v_j \tbeat v_k$. Then $\{v_i, v_j, v_k\}$ is agreeable, contradiction!
\end{nestedproof}

\begin{claim}
At both the \two\ triples $\{v_1, v_2, v_3\}$ and $\{v_4, v_5, v_6\}$ there are at least two vertices connected by a \one\ arrow to a vertex in the other \two\ triple.
\label{twovertices}
\end{claim}

\begin{nestedproof}
Suppose this would fail for $\{v_1, v_2, v_3\}$. Then two vertices, $v_i$ and $v_j$ with $\{i, j\} \in [4\setminus 1]^2$ would be connected to $\{v_4, v_5, v_6\}$ purely by \zero\ arrows. But there are six such connections to be made and since at most one \zero\ arrow may leave a vertex only the usage of five \zero\ arrows is allowed here---contradiction. The proof for the other case works totally analogously.
\end{nestedproof}

Now we are moving ever closer to a contradiction. By claim \ref{purelyone} we know that the \two\ triple $\{v_7, v_8, v_9\}$ is connected purely by \one\ arrows to one of the other two \two\ triples. Suppose $\{v_1, v_2, v_3\}$ is this triple. Since at most one \zero\ arrow may leave the \two\ triple $\{v_7, v_8, v_9\}$ in the direction of any vertex $v_i$ with $i \in 10\setminus 1$ and at most one \zero\ arrow can leave each vertex $v_i$ with $i \in 7\setminus 4$ we know that five arrows between $\{v_4, v_5, v_6\}$ and $\{v_7, v_8, v_9\}$ have to be \one. So in particular at least two vertices in $\{v_4, v_5, v_6\}$ have to have \one\ connections to vertices in the \two\ triple $\{v_7, v_8, v_9\}$. Using claim \ref{twovertices} and a variant of the pigeonhole principle we may choose a vertex $v_i$ with $i \in 7\setminus 4$ which has a \one\ connection both to a vertex $v_j$ with $j \in 4\setminus 1$ and to a vertex $v_k$ with $k \in 10\setminus 7$. Since we supposed the truth of claim \ref{purelyone} to be witnessed by $\{v_1, v_2, v_3\}$ we know that also $v_j$ and $v_k$ are connected by \one\ arrows. If we now remember that all \one\ arrows between $\{v_1, v_2, v_3\}$ and $\{v_7, v_8, v_9\}$ have to point from the latter to the former we know that in order for $\{v_i, v_j, v_k\}$ to be disagreeable the only possibility left is $v_i \obeat v_k \obeat v_j \obeat v_i$. But then we have $v_j \obeat v_0 \obeat v_i$ and $v_j \obeat v_i$ hence $\{v_0, v_i, v_j\}$ is agreeable. Contradiction. If the truth of claim \ref{purelyone} is witnessed by $\{v_4, v_5, v_6\}$ instead the proof works totally analogously. So this proves the theorem.
\end{aufz}
\end{proof}

\section{Strong agreeability}

\begin{definition}
A triple is called \emph{strongly agreeable} if and only if it is agreeable and does not contain any \two\ arrow. So it is strongly agreeable precisely if it is one of these:

\begin{tikzpicture}[
zero/.style={ultra thick, solid, ZERO},
one/.style={ultra thick, solid, ONE},
two/.style={ultra thick, solid, TWO}
]

\pgftransformscale{0.019}

\path
node (up) at (-339, 71) {}
node (left) at (-380,0) {}
node (right) at (-298,0) {};
\draw [->] (up) [one] -- (left);
\draw [->] (left) [one] -- (right);
\draw [->] (up) [one] -- (right);

\path
node (up) at (-226, 71) {}
node (left) at (-267,0) {}
node (right) at (-185,0) {};
\draw [->] (up) [zero] -- (left);
\draw [->] (left) [zero] -- (right);
\draw [->] (up) [zero] -- (right);

\path
node (up) at (-113, 71) {}
node (left) at (-154,0) {}
node (right) at (-72,0) {};
\draw [->] (up) [one] -- (left);
\draw [->] (left) [zero] -- (right);
\draw [->] (up) [one] -- (right);

\path
node (up) at (0, 71) {}
node (left) at (-41,0) {}
node (right) at (41,0) {};
\draw [->] (up) [zero] -- (left);
\draw [->] (left) [one] -- (right);
\draw [->] (up) [one] -- (right);

\path
node (up) at (113, 71) {}
node (left) at (72,0) {}
node (right) at (154,0) {};
\draw [->] (up) [zero] -- (left);
\draw [->] (left) [one] -- (right);
\draw [->] (up) [zero] -- (right);

\path
node (up) at (226, 71) {}
node (left) at (185,0) {}
node (right) at (267,0) {};
\draw [->] (up) [one] -- (left);
\draw [->] (left) [zero] -- (right);
\draw [->] (up) [zero] -- (right);
\end{tikzpicture}
\end{definition}

The proof of the following theorem is heavily inspired by \cite{974B0}.

\begin{theorem}
\label{kappalambdatheorem}
Let $\kappa$ be weakly compact and $\lambda\in\kappa\setminus\omega$ be a cardinal. The partition relation $\kappa\lambda l\arrows(\kappa\lambda m,n)$ holds true if and only if every edge-coloured digraph $C = \langle l, A, c\rangle$ with $\ran(c) = 2$ contains an independent set of size $m$ or there is a
subtournament $S$ of $C$ induced by a set of $n$ vertices such that all triples in $S$
are strongly agreeable.
\end{theorem}

\begin{proof}
\begin{aufz}
\item Towards a contradiction let us assume that the finite combinatorial characterization above holds true yet there is neither a $0$-homogeneous set of size $\kappa\lambda m$ nor a $1$-homogeneous $n$-tuple. Let furthermore $\chi:[\kappa\lambda l]^2\longrightarrow 2$ be any colouring. W.l.o.g. we may assume that $\pwi{\chi}{(\kappa(\lambda k+\alpha+1)\setminus\kappa(\lambda k+\alpha))}=1$ for any $\alpha<\lambda$ and $k<l$. We define a new colouring $\xi$ as follows:

\begin{align*}
\xi : [\kappa]^2&\longrightarrow \soafft{[\lambda l]^2}{4}\\
\{\alpha,\beta\}_<&\longmapsto\langle\{\gamma,\delta\}_<,2\chi(\{\kappa\gamma+\alpha,\kappa\delta+\beta\})+\chi(\{\kappa\gamma+\beta,\kappa\delta+\alpha\})\rangle\\
\end{align*}

Now we use fact \ref{partitionequivalence} in the form $\kappa \arrows (\kappa)^2_{2^\lambda}$ thereby finding an $X \in [\kappa]^\kappa$ homogeneous for $\xi$. Set $\psi:=\bigcup\pwi{\xi}{X}$. We have $\psi : [\lambda l]^2\longrightarrow 4$.

Let $\langle\langle b_j, c_j\rangle | j<2l(l-1)\rangle$ be an enumeration of $l^2\setminus\Delta$.

We are going to inductively construct for every $i<l$ a sequence of sets $\langle Y_i^j | j\leqslant 2l(l-1)\rangle$ by letting $Y_i^0:=\lambda(i+1)\setminus\lambda i$ and for any $j<l(l-1)$ letting $Y_{b_j}^{2j+1}:=B$ and $Y_{c_j}^{2j+1}:=C$ if there are sets $B\in[Y_{b_j}^{2j}]^\lambda$, $C\in [Y_{c_j}^{2j}]^\lambda$ with $\card{\soast{y\in C}{\psi(\{x,y\})\equiv 1(2)}}<\lambda$ for all $x\in B$. If there are no such sets let $Y_{b_j}^{2j+1}:=Y_{b_j}^{2j}$ and $Y_{c_j}^{2j+1}:=Y_{c_j}^{2j}$. Similarly let $Y_{b_j}^{2j+2}:=B$ and $Y_{c_j}^{2j+2}:=C$ if there are sets $B\in[Y_{b_j}^{2j+1}]^\lambda$, $C\in [Y_{c_j}^{2j+1}]^\lambda$ with $\card{\soast{y\in C}{\psi(\{x,y\})\in 4\setminus 2}}<\lambda$ for all $x\in B$. If there are no such sets let $Y_{b_j}^{2j+2}:=Y_{b_j}^{2j+1}$ and $Y_{c_j}^{2j+2}:=Y_{c_j}^{2j+1}$. Generally let $Y_i^{2j+2}:=Y_i^{2j}$ for $j<l(l-1)$ and $i\in l\setminus\{b_j,c_j\}$. Finally we let $Z_i:=Y_i^{2l(l-1)}$ for all $i<l$.

Now we define a digraph $D=\langle l, A, c\rangle$ by setting for $i,j<l$:

\begin{align}
\label{freedom0}
i\zbeat j\Longrightarrow&\exists\alpha\in Z_i:\card{\soast{\beta\in Z_j}{\psi(\{\alpha,\beta\})\equiv 1(2)}}=\lambda\\
\label{freedom1}
\Longleftrightarrow&\forall C\in[Z_j]^\lambda:\card{\soast{\eta\in Z_i}{\card{\soast{\nu\in C}{\psi(\{\eta,\nu\})\equiv 1(2)}}<\lambda}}<\lambda\\
\label{freedom2}
i\obeat j\Longrightarrow&\exists\alpha\in Z_i:\card{\soast{\beta\in Z_j}{\psi(\{\alpha,\beta\})\in 4\setminus 2}}=\lambda\\
\label{freedom3}
\Longleftrightarrow&\forall C\in[Z_j]^\lambda\card{\soast{\eta\in Z_i}{\card{\soast{\nu\in C}{\psi(\{\eta,\nu\})\in 4\setminus 2}}<\lambda}}<\lambda
\end{align}

We have $\Rightarrow$ instead of $\Leftrightarrow$ here for a reason. This definition neither necessarily tells us the direction of an arrow---since $i$ and $j$ could change roles---nor does it determine its colour uniquely---since we have $3\in\soast{x\in 4\setminus 2}{n\equiv 1(2)}$. We only have to follow the principle that whenever these conditions allow us to set an arrow we do so. The equivalences hold because of the definition of the $Z_i$'s. We will use the negations of (\ref{freedom0}) and (\ref{freedom2}) in the case where we define the $0$-homogeneous set and (\ref{freedom1}) and (\ref{freedom3}) in the case in which we define the $1$-homogeneous set.

After having so defined the digraph we may employ our finitary hypothesis and distinguish two cases:
\begin{aufz}
\item There is $N\in [l]^n$ such that all elements of $[N]^3$ are strongly agreeable. We are going to inductively construct a $1$-homogeneous set $H$. Since no strongly agreeable triple is cyclic we may choose an enumeration $e:n\longleftrightarrow N$ in such a way that $i<j$ if and only if $e(i)\beat e(j)$. Now we proceed by induction. We may, using (\ref{freedom1}) and (\ref{freedom3}), choose an $\alpha_i\in Z_{e(i)}$ such that
\begin{aufz}
\item $\forall j<i(e(i)\zbeat e(j)\rightarrow\psi(\{\alpha_j,\alpha_i)\equiv 1(2))$.

\item $\forall j<i(e(i)\obeat e(j)\rightarrow\psi(\{\alpha_j,\alpha_i)\in 4\setminus 2)$.

\item $\forall j\in n\setminus(i+1)(e(i)\zbeat e(j)\rightarrow\card{\soast{\beta\in Z_{e(j)}}{\psi(\{\alpha_0,\beta\})\equiv 1(2)}}=\lambda)$

\item $\forall j\in n\setminus(i+1)(e(i)\obeat e(j)\rightarrow\card{\soast{\beta\in Z_{e(j)}}{\psi(\{\alpha_0,\beta\})\in 4\setminus 2}}=\lambda)$
\end{aufz}

Now choose this end let $S\in [\kappa]^{2^n}$ be such that $\forall\eta\in S\forall\nu\in\eta\cap S\exists\zeta\in X\cap\eta\setminus\nu$. Let $\langle s_\eta|\eta<2^n\rangle$ be the increasing enumeration of $S$ and let $e:n\longleftrightarrow N$. Now let $k_{-2}:=0$, $k_{-1}:=2^n$\ox\ and $k_0:=2^{n-1}$. By induction we choose $\gamma_i\in X\cap s_{k_i}\setminus s_{k_i-1}$ and $k_{i+1}$ such that the following conditions are met.
\begin{aufz}
\item If $e(j)\zbeat e(i+1)$ then $k_{i+1}<k_j$.
\item If $e(j)\obeat e(i+1)$ then $k_{i+1}>k_j$.
\item $|k_{i+1}-k_j|\geqslant 2^{n-i-1}$ for every $j\leqslant i$.
\end{aufz}

This we do for every $i<n$. Note that we are always able to choose the parameters in this way because all elements of $[N]^3$ are strongly agreeable. Finally we may define $H:=\soast{\kappa\alpha_i+\gamma_i}{i<n}$. $H$ is $1$-homogeneous.

\item There is an independent $M\in [l]^m$. Again we have to distinguish two cases:
\begin{aufz}
\item $\lambda$ is regular. This case is easy. Let $Y_0:=\emptyset$ and for limit ordinals $\nu$ let $Y_\nu:=\bigcup_{\eta<\nu}X_\eta$. For any $\eta<\lambda$, $i\in M$ and an enumeration $f:m\longleftrightarrow M$ let

\begin{align}
\label{indcon}
Y_{m\eta+i+1}:=Y_{m\eta+i}\cup\{\min(\soast{\nu\in Z_{f(i)}}{\forall\zeta\in Y_{m\eta+i}:\psi(\{\nu,\zeta\})=0})\}
\end{align}

At the end let $Y:=Y_\lambda$. $M$ is independent so the formulas on the right of (\ref{freedom0}) and (\ref{freedom2}) are jointly negated.

The regularity of $\lambda$ now ensures that the minimum referred to in (\ref{indcon}) always exists.

\item $\lambda$ is singular. Here we have to be a bit more careful. We have to make sure that we do not spoil our chances of defining a homogeneous set before our inductive construction is finished. So we thin out the $Z_i$'s further as follows:

For each $i\in M$ let $\eta_i\in\lambda\cf(\lambda)+1\setminus\lambda$ be the least ordinal such that there exists an enumeration $e_i:\eta_i\longleftrightarrow Z_i$ such that

\begin{align*}
\nu\longmapsto\max(\soast{\card{\soast{\zeta\in Z_j}{\psi(\{e_i(\nu),\zeta\})>0}}}{j\in M\setminus\{i\}})
\end{align*}

is nondecreasing. Fix such $\eta_i$ and $e_i$. Similar to the case before we set $Y_0:=\emptyset$, for limit ordinals $\nu$ we let $Y_\nu:=\bigcup_{\zeta<\nu}X_\zeta$ and for $\nu<\lambda$, $i\in M$\ox\ and an enumeration $f:m\longleftrightarrow M$ we define

\begin{align}
\label{uglyformula}
Y_{m\nu+i+1}:=Y_{m\nu+i}\cup\{e_{f(i)}(\min(\soast{\zeta<\lambda}{\forall\rho\in Y_{m\nu+i}:\psi(\{e_{f(i)}(\zeta),\rho\})=0}))\}.
\end{align}

Once more set $Y:=Y_\lambda$. Again the independence of $M$ implies that the right sides of (\ref{freedom0}) and (\ref{freedom2}) are jointly negated. By rearranging the $Z_i$'s using the $e_i$'s we ensured that the minimum in (\ref{uglyformula}) always exists.
\end{aufz}

Let $\langle X_\eta|\eta<\lambda l\rangle$ be a partition of $X$ into $\lambda$ sets of size $\kappa$. We set $H:=\soast{\kappa\alpha+\gamma}{\alpha\in Y\wedge\gamma\in X_\alpha}$. A little contemplation shows that $H$ is $0$-homogeneous for $\chi$.
\end{aufz}

\item Suppose that the finite combinatorial characterization above fails, i.e. there is an edge-$2$-coloured digraph $D=\langle l, A, c\rangle$ neither containing an independent set of size $m$ nor a set of size $n$ only inducing strongly agreeable triples. We define a colouring $\chi:[\kappa\lambda l]^2\longrightarrow 2$ as follows.

\begin{align*}
\chi:[\kappa\lambda l]^2&\longrightarrow 2\\
\{\kappa(\lambda j+\alpha) + \gamma,\kappa(\lambda k+\beta)+\delta\}_<&\longmapsto\begin{cases}
1\text{ if }j\zbeat k\wedge\alpha<\beta\wedge\delta<\gamma\\
1\text{ if }k\zbeat j\wedge\beta<\alpha\wedge\gamma<\delta\\
1\text{ if }j\obeat k\wedge\alpha<\beta\wedge\gamma<\delta\\
1\text{ if }k\obeat j\wedge\beta<\alpha\wedge\delta<\gamma\\
0\text{ else.}
\end{cases}
\end{align*}
\begin{aufz}
\item Suppose that $\chi$ would admit a $0$-homogeneous set $X\in[\kappa\lambda l]^{\kappa\lambda m}$. Then there would have to be a set $Y\in[l]^m$ such that $\otyp(X\cap\kappa\lambda(k+1)\setminus\kappa\lambda k) = \kappa\lambda$ for any $k\in Y$. Since $D$ contains no independent set of size $m$ we may choose $\{j,k\}\in[Y]^2$ with $j\beat k$. Let $\alpha:=\min\soast{\eta<\lambda}{\otyp(X\cap\kappa(\lambda j+\eta+1)\setminus\kappa(\lambda j+\eta))=\kappa}$ and let $\beta:=\min\soast{\eta\in\lambda\setminus(\alpha+1)}{\otyp(X\cap\kappa(\lambda k+\eta+1)\setminus\kappa(\lambda k+\eta)=\kappa}$.

We distinguish two cases:
\begin{aufz}
\item $c(\{j,k\})=0$. Let $\gamma:=\min\soast{\eta<\kappa}{\kappa(\lambda k+\beta)+\eta\in X}$ and let $\delta:=\min\soast{\eta\in\kappa\setminus(\gamma+1)}{\kappa(\lambda j+\alpha)+\eta\in X}$. Obviously $\chi(\{\kappa(\lambda j+\alpha)+\delta,\kappa(\lambda k+\beta)+\gamma\})=1$ contradicting the $0$-homogeneity of $X$.

\item $c(\{j,k\})=1$. Let $\gamma:=\min\soast{\eta<\kappa}{\kappa(\lambda j+\alpha)+\eta\in X}$ and let $\delta:=\min\soast{\eta\in\kappa\setminus(\gamma+1)}{\kappa(\lambda k+\beta)+\eta\in X}$. Obviously $\chi(\{\kappa(\lambda j+\alpha)+\gamma,\kappa(\lambda k+\beta)+\delta\})=1$ contradicting the $0$-homogeneity of $X$.
\end{aufz}

\item Now suppose that $\chi$ would admit a $1$-homogeneous $n$-tuple $X\in[\kappa\lambda l]^n$. Per definitionem of $\chi$ we have that $\{\kappa(\lambda j+\alpha)+\gamma,\kappa(\lambda k+\beta)+\delta\}\in [X]^2$ for $j,k<l$; $\alpha,\beta<\lambda$ and $\gamma,\delta<\kappa$ implies $j\ne k$; $\alpha\ne\beta$ and $\gamma\ne\delta$. Consider $Y:=\soast{k<l}{\exists \eta<\kappa\lambda:\kappa\lambda+\eta\in X}$. We have $Y\in[l]^n$ and in fact $Y$ induces a subdigraph $S$ of $D$. W.l.o.g. we may assume that $S$ is a tournament because if $\{j,k\}\in[Y]^2$ is independent we have that $X\cap((\kappa\lambda(j+1)\setminus\kappa\lambda j)\cup(\kappa\lambda(k+1)\setminus\kappa\lambda k)\subset X$ is $0$-homogeneous.

So $Y$ is a tournament and contains a triple $\{i,j,k\}\in[Y]^3$ which is not strongly agreeable. Per definitionem of $Y$ we may conclude that there are $\alpha,\beta,\gamma<\lambda$ and $\delta,\eta,\zeta<\kappa$ such that $\kappa(\lambda i+\alpha)+\delta,\kappa(\lambda j+\beta)+\eta,\kappa(\lambda k+\gamma)+\zeta\in X$. We may distinguish three cases:
\begin{aufz}
\item It is a $3$-cycle. Assume $i\beat j\beat k\beat i$. Since $X$ is $1$-homogeneous we may, per definitionem of $\chi$, conclude that $\alpha<\beta$, $\beta<\gamma$\ox\ and $\gamma<\alpha$. \con

\item We have $i\zbeat j$, $j\zbeat k$\ox\ and $i\obeat k$. By $1$-homogeneity of $X$ and per definitionem of $\chi$ we may conclude that $\delta>\eta$, $\eta>\zeta$\ox\ and $\delta<\zeta$. \con

\item We have $i\obeat j$, $j\obeat k$\ox\ and $i\zbeat k$. Again we may use the $1$-homogeneity of $X$ and the definition of $\chi$ to conclude that $\delta<\eta$, $\eta<\zeta$\ox\ and $\delta>\zeta$. \con
\end{aufz}

Note that modulo a renaming of $i$, $j$\ox\ and $k$ these three cases exhaust all possibilities.
\end{aufz}

This finishes the second part of the proof.
\end{aufz}
\end{proof}

Similar to before we denote by $S_n$ the class of graphs $G$ on $n$ vertices such that all triples in $G$ are strongly agreeable and by $r(I_m, S_n)$ the smallest natural number $l$ such that all digraphs on $l$ vertices either contain an independent set of size $m$ or an induced subgraph on $n$ vertices all triples of which are strongly agreeable. As before we may use this notation to condense the equivalences proved above into a single formula:

\begin{align*}
r(\kappa\lambda m, n) = \kappa\lambda r(I_m, S_n) \text{ for } \kappa \text{ weakly compact and any cardinal } \lambda\in\kappa\setminus\omega.
\end{align*}

Analogous to \cite{989B0} we can prove another nice theorem.
\begin{theorem}
\label{MAtheorem}
If $\MA_{\aleph_1}$ holds true then $r(\omega_1\omega m,n)=\omega_1\omega r(I_m,S_n)$.
\end{theorem}

Let $\chi : [\omega_1\omega l]^2\longrightarrow 2$ be any colouring. We define $E := \soast{P \in [\omega_1\omega l]^2}{\chi(P)=1}$.

Let $Z := \bigcup_{\alpha < \On} [\omega_1(\alpha+1)\setminus\omega_1\alpha]^{\aleph_1}$

For $X \in Z$ let $\alpha_X < \On$ be such that $X \in [\omega_1(\alpha_X+1)\setminus\omega_1\alpha_X]^{\aleph_1}$. Let $A_X := \soast{\xi<\omega_1}{\omega_1\alpha_X+\xi \in X}$. We say that the pair $\langle X, Y\rangle$ is
\begin{aufz}
\item in constellation $\zb$ iff $X, Y \in Z$ and\\

$\forall B \in [X]^{\aleph_1}, C \in [Y]^{\aleph_1} \exists \gamma\in A_C, \beta \in A_B\setminus (\gamma + 1) : \{\omega_1\alpha_X+\beta, \omega_1\alpha_Y+\gamma\} \in E$.

\item in constellation $\ob$ iff $X, Y \in Z$ and\\

$\forall B \in [X]^{\aleph_1}, C \in [Y]^{\aleph_1} \exists \beta\in A_B, \gamma \in A_C\setminus (\beta + 1) : \{\omega_1\alpha_X+\beta, \omega_1\alpha_Y+\gamma\} \in E$.

\item \emph{free} iff $\forall B \in [X]^{\aleph_1}, C \in [Y]^{\aleph_1} \exists D \in [B]^{\aleph_1}, F \in [C]^{\aleph_1} \forall \delta \in D, \eta \in F : \{\delta, \eta\} \notin E$.
\end{aufz}

Let us call a set $S \subset \On$ weakly independent if for all $\alpha, \beta < \On$ with $X = S \cap \omega_1(\alpha + 1)\setminus\omega_1\alpha$ and $Y = S \cap \omega_1(\beta + 1) \setminus \omega_1\beta$ and $\card{X}=\card{Y}=\aleph_1$ we have that $\langle X, Y\rangle$ is free.

For the proof we are going to use theorem 3.1 from \cite{989B0}. By limiting its scope to the case $\alpha = \omega$ we get the following theorem.

\begin{theorem}
\label{JEBsTheorem}
Assume $\MA_{\aleph_1}$. Let $E \subset [\omega_1\omega]^2$ be a graph and let $W \in [\omega_1\omega]^{\omega_1\omega}$ be weakly independent. Then there exists $H \in [W]^{\omega_1\omega}$ which is independent such that $\forall n < \omega : \card{H \cap \omega_1(n+1)\setminus \omega_1 n}\in\{0,\aleph_1\}$ and $H \cap \omega_1(n+1)\setminus \omega_1 n\supsetneq\emptyset$ implies $\card{W \cap \omega_1(n+1)\setminus \omega_1 n} \leqslant \aleph_0$ for all natural $n$.
\end{theorem}
Now we prove theorem \ref{MAtheorem}.

\begin{proof}
Even without $\MA_{\aleph_1}$ it is easy to see that $r(\omega_1\omega m, n) \geqslant \omega_1\omega r(I_m, S_n)$. Let $D$ be a $2$-coloured digraph showing $l < r(I_m, S_n)$. Now consider the following colouring:

\begin{align*}
\chi: [\omega_1\omega l]^2 & \longrightarrow 2\\
\{\omega_1\omega h + \omega_1 j + \alpha, \omega_1\omega i + \omega_1 k + \beta\}_<& \longmapsto \begin{cases}
1 \text{ if } h \zbeat i \wedge j < k \wedge \beta < \alpha\\
1 \text { if } i \zbeat h \wedge k < j \wedge \alpha < \beta\\
1 \text{ if } h \obeat i \wedge j < k \wedge \alpha < \beta\\
1 \text { if } i \obeat h \wedge k < j \wedge \beta < \alpha\\
0 \text{ else.}
\end{cases}
\end{align*}

Analogously to the proof of theorem \ref{kappalambdatheorem} this now shows $\omega_1\omega l \doesntarrow (\omega_1\omega m, n)$.

Now assume $\MA_{\aleph_1}$. Let $l := r(I_m, S_n)$. We are going to show $\omega_1\omega l \arrows (\omega_1\omega m, n)$ by adapting the method of \cite{989B0} for our purposes.

Note that if a pair $\langle X, Y\rangle$ is not free there are $X' \in [X]^{\aleph_1}$ and $Y' \in [Y]^{\aleph_1}$ such that $\langle X', Y'\rangle$ is in one of the two other constellations. If at least one of these three properties applies to a pair then we call it \emph{decided}.

We are choosing $l$-tuples $\langle X^F_0, \dots, X_{l-1}^F\rangle$ with $F \in [\omega]^{<\omega}$ by recursion on $\max(F)$. \begin{aufz}

\item If $\card{F} = 1$ and $i < l$ let $X^F_i \in [\omega_1(\omega i + \bigcup F + 1)\setminus\omega_1(\omega i + \bigcup F)]^{\aleph_1}$ such that for all $G \in [\bigcup F]^{<\omega}$ and $j \in l\setminus\{i\}$ either $\langle X^G_j, X^F_i\rangle$ is free or there is an $X \in [X^G_j]^{\aleph_1}$ such that $\langle X, X^F_i\rangle$ is in one of the other two constellations above. This can be achieved by thinning out repeatedly.

\item If $\card{F} \in \omega \setminus 2$ let $\mu := \max(F)$ and $G := F \setminus \{\mu\}$. Now by again thinning out repeatedly for all $i < l$ choose $X^F_i\subset X^G_i$ such for all $j\in l\setminus\{i\}$ we have that $\langle X^F_i, X^{\{\mu\}}_j\rangle$ is decided.
\end{aufz}

\begin{claim}
\label{techclaim}
For every natural $f$ there is an $F \in [\omega]^f$ and an $H \in [\omega]^\omega$ such that for all $g, h < l$ the facts both whether or not for $x, y \in [H]^i$ with $i \in F$ and $\max(x) < \min(y)$ we have that $\langle X^x_g, X^y_h\rangle$ is free or not free and what the constellation of a pair $\langle X_g^x, X^y_h\rangle$ with $\min(y)\in x$ and $x, y \in [H]^i$ is only depends on $g$ and $h$.
\end{claim}

\begin{nestedproof}
We repeatedly use Ramsey's theorem to arrive at $F$ and $H$.

Using Ramsey's theorem we thin $\omega$ to some $H_0 \in [\omega]^\omega$ on which the freeness of some such $\langle X^x_g, Y^y_h\rangle$ for $x, y \in [H_0]^1$ only depends on $\{g, h\}$ and whether $\max(x) < \min(y)$ or $\max(y) < \min(x)$.

Because of the nature of our recursive construction the constellation of a pair $\langle X_g^x, X^y_h\rangle$ with $\min(y) \in x$ is determined by the constellation of the pair $\langle X_g^{x\cap(\min(y)+1)}, X^{\{\min(y)\}}_h\rangle$. The constellation of this pair however can be made canonical as well using Ramsey's theorem. So inductively do this for $x\cap\min(y)=i$ to arrive at some $H_i'$ for $i<\omega$ and after that find an $H_{i+1} \in [H_i']^\omega$ on which the freeness/non-freeness for $x, y \in [H_{i+1}]^{i+1}$ with $\max(x)<\min(y)$ is canonical, i.e. only depends on $g$ and $h$. This yields a sequence $H_0 \supset H_0' \supset H_1 \supset H_1' \supset\dots$.

Since there are only finitely many possible patterns of both freeness/non-freeness and constellations we may find an $F \in [\omega]^f$ and an $H$ in the sequence which is sufficiently thin such that these patterns are the same for all $i \in F$.
\end{nestedproof}

Use the claim to find $F$ of a sufficiently( for what folllows) large finite size $2L$, say $2^{n+1}-2$.

Now we define a digraph by letting there be arrows as follows:

\begin{align*}
g\zbeat h \text{ implies } \langle X^x_g, X^{\{\max(x)\}}_h\rangle \text{ generally is in constellation } \zb,\\
g\obeat h \text{ implies }\langle X^x_g, X^{\{\max(x)\}}_h\rangle \text{ generally is in constellation } \ob.
\end{align*}

Moreover whenever one of these conditions applies an arrow is set.

Now first suppose there were an independent set $I \in [l]^m$. Let $\langle h_i | i < \omega\rangle$ be the ascending enumeration of $H$ and pick some $f\in F$. Let $Z := \soast{\{h_{fi}, \dots, h_{f(i+1)-1}\}}{i < \omega}$. Then the union of all the $X^x_g$'s for $g\in I$ and $x \in Z$ is a weakly independent set of size $\omega_1\omega m$. This is because by our recursive construction and by definition of our digraph $\langle X^{x\setminus\{\max(x)\}}_g, X^y_h\rangle$ is free whenever $\max(x) \in y$ and there is no arrow between $g$ and $h$ . We now may use theorem \ref{JEBsTheorem} to thin this out to some truly independent set of the same order-type, i.e. some set of order-type $\omega_1\omega m$ which is homogeneous for $\chi$ in colour $0$.

Now suppose that our digraph contains no independent set of size $m$. Since we chose $l$ to be $r(I_m, S_n)$ we know that there is a strongly agreeable set $S \in [l]^n$.

Let $\langle f_i | i < L\rangle$ be an ascending enumeration of $F$. Define $a_i := \{h_{f_i}, \dots, h_{f_L}\rangle$ and let $A^i_g := X^{a_i}_g$.

\begin{claim}
Let $g, h \in S$, $\{i, j\} \in [L]^2$\ox\ and $i < j$. Now if $g \zbeat h$ then $\langle A^i_g, A^j_h\rangle$ is in constellation \zb\ and if $g \obeat h$ then it is in constellation \ob.
\end{claim}

\begin{nestedproof}
Consider $x := \{h_{f_i}, \dots, h_{f_{j-1}}\}$ and $y := \{h_{f_j}, \dots, h_{f_{2j-i-1}}\}$.

We know that $\langle X^x_g, X^y_h\rangle$ is not free. Because of $X^y_g \subset X^{\{h_{f_j}\}}_h$ we have that $\langle X^x_g, X^{\{h_{f_j}\}}\rangle$ is not free. By our recursive construction we get that $\langle X^{x \cup \{h_{f_j}\}}_g, X^{\{h_{f_j}\}}_h\rangle$ is decided. The colour of the arrow between $g$ and $h$ in $S$ tells us how. Since $A^i_h \subset X^{x \cup \{h_{f_j}\}}_h$ and $A^j_g \subset X^{\{h_{f_j}\}}_g$ the claim follows.
\end{nestedproof}

Now we are going to construct the complete graph $G \in [\omega_1\omega l]^n$ with $[G]^2 \subset E$. To this end let $\langle s_i | i < n\rangle$ be an enumeration of $S$ such that $s_i \obeat s_j$ or $s_j \zbeat s_i$ for all $i < j < n$. This is possible since all triples in $S$ are strongly agreeable. To start the construction now, i.e. in step $0$ choose the middle element $j_0$ of $L$. This divides $L$ into two intervals of size, say, $2^{n-1}-1$. Let us call $\{-1, L\}$ \emph{the border}.

Now inductively in step $i < n$ we haven chosen $j_k$'s with $k < i$ such that the minimal distance between either of them and between one of them and the border is $2^{n-i}-1$. Now we choose an $j_{i+1}$ such that the minimal distance between two of the $j_{k+1}$'s for $k \leqslant i$ and between of them and the border is $2^{n-i-1}-1$ and such that $j_{i+1}>j_k$ iff $s_k\beat s_i$ for all $k \leqslant i$. Again strong agreeability allows us to do this.

Finally we choose in step $0$ an ordinal $\alpha_0 \in A^{j_0}_{s_0}$ such that for all $i \in n\setminus 1$ we have $\card{\soast{\alpha \in A^{j_i}_{s_i}}{\{\alpha_0, \alpha\} \in E}} = \aleph_1$. Such an ordinal exists for otherwise there would be an $i < n$ and uncountably many $\alpha\in A^{j_0}_{s_0}$ which are only connected to countably many things in $A^{j_i}_{s_i}$ by $E$. Because of the way in which we sorted the elements of $S$ this would contradict the constellation of $\langle A^{j_0}_{s_0}, A^{j_i}_{s_i}\rangle$. Now thin out the $A^{j_i}_{s_i}$ for $i \in n\setminus 1$ to those uncountable subsets of ordinals the unordered pairs with $\alpha_0$ of which are in $E$. Now generally in step $i < n$ choose an ordinal $\alpha_i \in A^{j_i}_{s_i}$ such that $\{\alpha_k, \alpha_i\} \in E$ for all $k < i$ and thin out the $A^{j_k}_{s_k}$'s for $k \in n\setminus (i+1)$. In this way we may construct our complete graph of size $n$, i.e. our set homogeneous for $\chi$ in colour $1$.
\end{proof}

Now we may take a look at the easiest examples.

\begin{theorem}
$r(I_2, S_3) = 6.$
\label{simplestnontrivialexample}
\end{theorem}

\begin{proof}
\mbox{}\\
\begin{minipage}{0.545\linewidth}
\begin{tikzpicture}[scale=2,
zero/.style={ultra thick, solid, ZERO},
one/.style={ultra thick, solid, ONE},
vertex/.style={ultra thick, circle, draw=dark grey, fill=very light grey, minimum size=3cc}]

\pgftransformscale{0.6}

\draw (0,0) arc (0:72:3) [white] node [vertex] (a) {}
arc (72:144:3) node (b) [vertex] {}
arc (144:216:3) node (c) [vertex] {}
arc (216:288:3) node (d) [vertex] {}
arc (288:0:3) node (e) [vertex] {};

\foreach \from/\to in {a/b,b/c,c/d,d/e,e/a,c/a,d/b,e/c,a/d,b/e}
\draw [->] (\from) [zero] -- (\to);

\foreach \from/\to in {a/d,b/e,c/a,d/b,e/c}
\draw [->] (\from) [one] -- (\to);
\end{tikzpicture}
\end{minipage}
\hfill
\begin{minipage}{0.364\linewidth}
The lower bound is provided by the counterexample on the right, for the upper bound one can argue as follows:

Consider any vertex $v$ in a six-element tournament---if the digraph is not a tournament we are finished immediately. Now $\zinh(v) + \zonh(v) + \oinh(v) + \oonh(v) \geqslant 5$ so we may conclude using the pigeonhole principle that one of these four neighborhoods contains at least two elements. But then $v$ together with this neighborhood is strongly agreeable.
\end{minipage}\\
\end{proof}

\begin{theorem}
$r(I_3, S_3) = 15$.
\end{theorem}

\begin{proof}
\hspace{-5cc}
\begin{tikzpicture}[scale=2,
zero/.style={ultra thick, solid, ZERO},
one/.style={ultra thick, solid, ONE},
vertex/.style={ultra thick, circle, draw=dark grey, fill=very light grey, minimum size=3cc}]

\pgftransformscale{0.95}

\draw (0,0) arc (0:25.7:3) [white] node [vertex] (a) {}
arc (25.7:51.4:3) node (b) [vertex] {}
arc (51.4:77.1:3) node (c) [vertex] {}
arc (77.1:102.9:3) node (d) [vertex] {}
arc (102.9:128.6:3) node (e) [vertex] {}
arc (128.6:154.3:3) node (f) [vertex] {}
arc (154.3:180:3) node (g) [vertex] {}
arc (180:205.7:3) node (h) [vertex] {}
arc (205.7:231.4:3) node (i) [vertex] {}
arc (231.4:257.1:3) node (j) [vertex] {}
arc (257.1:282.9:3) node (k) [vertex] {}
arc (282.9:308.6:3) node (l) [vertex] {}
arc (308.6:334.3:3) node (m) [vertex] {}
arc (334.3:0:3) node (n) [vertex] {};

\foreach \from/\to in {a/b,b/c,c/d,d/e,e/f,f/g,g/h,h/i,i/j,j/k,k/l,l/m,m/n,n/a}
\draw [->] (\from) [zero] -- (\to);

\foreach\from/\to in {e/a,f/b,g/c,h/d,i/e,j/f,k/g,l/h,m/i,n/j,a/k,b/l,c/m,d/n}
\path (\from) edge [->, bend right=20, zero] (\to);

\foreach \from/\to in {a/c,b/d,c/e,d/f,e/g,f/h,g/i,h/j,i/k,j/l,k/m,l/n,m/a,n/b}
\path (\from) edge [->, bend right=60, one] (\to);

\foreach\from/\to in {d/a,e/b,f/c,g/d,h/e,i/f,j/g,k/h,l/i,m/j,n/k,a/l,b/m,c/n}
\draw [->] (\from) [one] -- (\to);

\end{tikzpicture}

The lower bound is given by the counterexample above. It can be viewed as a digraph on $\mathbb{Z}_{14}$ where $v\zbeat v+1$, $v\obeat v+2$, $v\obeat v-3$\ox\ and $v\zbeat v-4$ for all $v\in\mathbb{Z}_{14}$.

\end{proof}

The upper bound is established by theorem \ref{secondboundforstronglyagreeables} below.

\begin{corollary}
\ \\Let $\kappa$ be the least weakly compact cardinal. Then $r(\kappa\aleph_{\varepsilon_{\omega 7 + 9}}3,3)=\kappa\aleph_{\varepsilon_{\omega 7 + 9}}15$.
\end{corollary}

\section{Upper bounds}

\subsection{Upper bounds for independent versus strongly agreeable sets.}

We now are going to provide several general upper bounds for various Ramsey numbers. We start with something easy---an upper bound for the Ramsey numbers $r(I_m, S_3)$.

\begin{lemma}
\label{inductivestep0}
$r(I_{m + 1}, S_3) \leqslant r(I_m, S_3) + 4m + 1$ for all $m \in \omega \setminus 2$.
\end{lemma}

\begin{proof}
Let $D$ be any edge-bicoloured digraph on $r(I_m, S_3) + 4m + 1$ vertices. Fix any vertex $v$. Since $\card{\zonh(v)}+\card{\zinh(v)}+\card{\oonh(v)}+\card{\oinh(v)}+\card{\ind(v)}=r(I_m, S_3) + 4m$ it follows that one of the following cases has to obtain:
\begin{aufz}
\item $\card{\zonh(v)}\geqslant m+1$.
\item $\card{\zinh(v)}\geqslant m+1$.
\item $\card{\oonh(v)}\geqslant m+1$.
\item $\card{\oinh(v)}\geqslant m+1$.
\item $\card{\ind(v)}\geqslant r(I_m, S_3)$.
\end{aufz}
In the first four cases either the respective neighbourhood is independent or two vertices in it are connected by an arrow in which case they form an $S_3$ together with $v$. In the last case the neighbourhood either contains an $S_3$ or an independent set of size $m$ which together with $v$ would form an independent set of size $m+1$.
\end{proof}

\begin{theorem}
For all $m \in\omega\setminus 2$ we have $r(I_m, S_3) \leqslant m(2m - 1)$.
\label{firstboundforstronglyagreeables}
\end{theorem}

\begin{proof}
The base case is established by theorem \ref{simplestnontrivialexample}. For the inductive step by lemma \ref{inductivestep0} we only have to check that
\begin{align*}
(m+1)(2(m+1)-1)=m(2m-1)+4m+1
\end{align*}
which is true.
\end{proof}

We continue now by providing an upper bound for $r(I_2, S_n$).

\begin{lemma}
For all $n\in\omega\setminus 3$ we have $r(I_2, S_{n+1}) \leqslant 4r(I_2, S_n) - 2$.
\end{lemma}

\begin{proof}
Let $T = \langle V, A, c\rangle$ be a bicoloured tournament on $4r(I_2, S_n) - 2$ vertices. Fix any vertex $v \in V$. Since $\card{\zinh(v) + \zonh(v) + \oinh(v) + \oonh(v)} = 4r(I_2, S_n) - 3$ by the pigeonhole principle one of these neighbourhoods $N$ has cardinality at least $r(I_2, S_n)$. Hence there is an $S \in [N]^n$ in which all triples are strongly agreeable. So if we set $X := S \cup \{v\}$ we have $X \in [V]^{n+1}$ and all triples in $X$ are strongly agreeable.
\end{proof}

\begin{theorem}
For any $n\in\omega\setminus 3$ we have
\begin{align*}
r(I_2, S_n) \leqslant \frac{4^{n-1} + 2}{3}.
\end{align*}
\label{secondboundforstronglyagreeables}
\end{theorem}

\begin{proof}
By induction. The base case is given by theorem \ref{simplestnontrivialexample}. The induction step also works:
\begin{align*}
4\cdot\frac{4^{n-1} + 2}{3} - 2 = \frac{4^{n-1+1} + 4\cdot 2 - 2\cdot 3}{3} = \frac{4^{n+1-1}+ 2(4-3)}{3}
\end{align*}
\end{proof}

Again, note that for any integer $n$ also $\frac{4^{n-1}+2}{3}$ is an integer. It is easy to see by induction that $4^{n-1}\equiv 1(3)$ for any $n\in\omega\setminus 1$. Hence $4^{n-1}+2\equiv 0(3)$ for any $n\in\omega\setminus 1$.

Now we want to give a general upper bound for $r(I_m, S_n)$.

\begin{lemma}
$r(I_{m+1}, S_{n+1})\leqslant r(I_m,S_{n+1})+4r(I_{m+1},S_n)-3$ for all $m \in \omega \setminus 2$ and all $n \in \omega \setminus 3$.
\label{mainlemma}
\end{lemma}

\begin{proof}
Immediate.
\end{proof}

\begin{theorem}For all $m\in\omega\setminus 2$ and all $n\in\omega\setminus 3$ we have $r(I_m, S_n) \leqslant u(m, n)$ where
\begin{align}
u(m, n):=&\frac{1}{4}\Big( 3+\sum_{i=0}^{n-1}\binom{i+m-2}{i}4^i\Big).
\end{align}
\end{theorem}

\begin{proof}
By induction. First we verify that $u(2, 3) = 6$. Recall that we have $r(I_2, S_3) = 6$ by theorem \ref{simplestnontrivialexample}. Next we check that for $u(2,n)$ agrees with the formula from theorem \ref{secondboundforstronglyagreeables}. We have

\begin{align*}
u(2,n)=\frac{1}{4}\big(3+\sum_{i=0}^{n-1} 4^i\big).
\end{align*}

Consider $v(n):=u(2,n+1)-u(2, n)$. It is easy to see that $v(n)=4^{n-1}$. But it is also easy to see that

\begin{align*}
\frac{4^n+2}{3}-\frac{4^{n-1}+2}{3}=4^{n-1}.
\end{align*}

Induction on $n$ provides what was demanded.

It is even easier to see, using well-known properties of the binomial coefficients that $u(m, 3)=m(2m-1)$.

Now we can check that in fact

\begin{align}
\label{upeq}
u(m+1,n+1)=u(m,n+1)+4u(m+1,n)-3.
\end{align}

Then with lemma \ref{mainlemma} the theorem will follow. The following calculation proves (\ref{upeq}).

\begin{align*}
&u(m,n+1)+4u(m+1,n)-3\\
=\yq&\frac{1}{4}\Big(3+\sum_{i=0}^n\binom{i+m-2}{i}4^i+4\sum_{i=0}^{n-1}\binom{i+m-1}{i}4^i\Big)\\
=\yq&\frac{1}{4}\Big(3+\sum_{i=0}^n\binom{i+m-2}{i}4^i+\sum_{i=1}^n\binom{i+m-2}{i-1}4^i\Big)\\
=\yq&\frac{1}{4}\Big(4+\sum_{i=1}^n\Big(\binom{i+m-2}{i}+\binom{i+m-2}{i-1}\Big)4^i\Big)\\
=\yq&\frac{1}{4}\Big(4+\sum_{i=1}^n\binom{i+m-1}{i}4^i\Big)\\
=\yq&\frac{1}{4}\Big(3+\sum_{i=0}^n\binom{i+m-1}{i}4^i\Big)\\
=\yq&u(m+1,n+1).
\end{align*}
\end{proof}

\subsection{Upper bounds for independent versus agreeable sets.}

We are now going to provide an upper bound for the Ramsey numbers $r(I_n, A_3)$.

\begin{lemma}
\label{inductivestep1}
$r(I_{n + 1}, A_3) \leqslant r(I_n, A_3) + 2r(I_{n + 1}, L_3) + 4n - 1$ for all $n \in \omega \setminus 2$.
\end{lemma}

\begin{proof}
Let there be an edge-tricoloured digraph $\langle V, A\rangle$ with\linebreak $r(I_n, A_3) + 2r(I_{n + 1}, L_3) + 4n - 1$ vertices. Pick any vertex $v$ such that $\zonh(v)$ is maximal. Note that $\card{V \setminus \{v\}} = r(I_n, A_3) + 2r(I_{n + 1}, L_3) + 4n - 2$! Now we may distinguish several cases. Note that necessarily at least one of them must hold!

\begin{aufz}

\item $\card{\zinh(v)} \geqslant n + 1$. Since we chose $v$ to have maximum \zero\ out-degree we might immediately hand over to the next case.

\item $\card{\zonh(v)} \geqslant n + 1$.
In this case either there are two vertices $w, x \in \zonh(v)$ which are connected by an arrow---in which case $\{v, w, x\}$ is an agreeable triple---or $\zonh(v)$ forms an independent set of size $n + 1$. In both cases we found a set homogeneous in the sought-after-sense.

\item $\card{\oinh(v)} \geqslant r(I_{n + 1}, L_3)$.
If there are $w, x \in \oinh(v)$ which are connected by either a \zero\ or a \one\ arrow then $\{v, w, x\}$ is an agreeable triple. So let us assume that all pairs of vertices in $\oinh(v)$ are either connected by \two\ arrows or not connected at all. But now our case condition implies that we either have an independent set of size $n + 1$ or a transitive \two\ triple---which is agreeable.

\item $\card{\oonh(v)} \geqslant r(I_{n + 1}, L_3)$.
Works exactly like the case before.

\item $\card{\tinh(v)} \geqslant n + 1$.
This and the next case work like the second. Either there are two vertices $w, x \in \tinh(v)$ which are connected by an arrow---then $\{v, w, x\}$ is agreeable---or $\tinh(v)$ is an independent set of size $n + 1$.

\item $\card{\tonh(v)} \geqslant n + 1$. Again either there are two vertices which are connected which then form an agreeable triple together with $v$ or $\tonh(v)$ is independent.

\item $\card{\ind(v)} \geqslant r(I_n, A_3)$.
Now in the last case $\ind(v)$ either contains an agreeable triple or is itself independent of size $n$. In the first case we are finished immediately and in the second $\ind(v) \cup \{v\}$ is an independent set of size $n + 1$.
\end{aufz}
\end{proof}

\begin{corollary}
$r(\omega^2(n + 1), 3) \leqslant r(\omega^2 n, 3) + \omega2r(\omega(n + 1), 3) + \omega^2(4n - 1)$ for all $n \in \omega \setminus 2$.
\end{corollary}

\begin{theorem}
\label{upperbound}
For all $n \in \omega \setminus 2$ we have
\begin{align}
\label{upperformula}
r(I_n, A_3) \leqslant \frac{(2n+1)(n^2 + 4n - 6)}{3}.
\end{align}
\end{theorem}

\begin{proof}
By induction. The basis case is $n = 2$, formula \ref{upperformula} is verified by theorem \ref{firstnontrivial} in it. So let us assume that the formula holds up to some $n \in \omega \setminus 2$. Using lemma \ref{inductivestep1} we get
\begin{align*}
& r(I_{n + 1}, A_3) \leqslant r(I_n, A_3) + 2r(I_{n + 1}, L_3) + 4n - 1.
\end{align*}
Larson and Mitchell---see Lemma 4.1 of \cite{997LM0}---showed that $r(I_n, L_3) \leqslant n^2$, so we may bound this from above as follows:
\begin{align*}
\leqslant & \frac{(2n + 1)(n^2 + 4n - 6)}{3} + 2(n + 1)^2 + 4n - 1\\
= & \frac{(2(n + 1) + 1)((n + 1)^2 + 4(n + 1) - 6)}{3}
\end{align*}
\end{proof}

Note that the fraction in theorem \ref{upperbound} is an integer if $n$ is one. If $n\equiv 1(3)$ then $2n+1\equiv 0(3)$ and we are fine. If, however, $n\equiv 0(3)$ then $n^2+4n-6\equiv 0(3)$ and if $n\equiv 2(3)$ then $n^2+4n-6\equiv 0(3)$ too.

\section{A semiconstructive lower bound}

\begin{theorem}
$r(I_{mn+1}, A_3) > (r(I_{m+1}, L_3)-1)(r(I_{n+1},L_3)-1)$.
\end{theorem}

\begin{proof}
Let $D_0 = \langle V_0, A_0\rangle$ be a digraph on $r(I_m, L_3)-1$ vertices neither containing a transitive triple nor an independent set of size $m+1$ and let $D_1 = \langle V_1, A_1\rangle$ be a digraph on $r(I_n, L_3)-1$ vertices neither containing a transitive triple nor an independent set of size $n+1$. Then consider the coloured digraph $D_0*D_1 := \langle V_0\times V_1, A_2, c\rangle$ where $A_3 = \soast{\langle\langle v_0, v_1\rangle, \langle v_2, v_1\rangle\rangle}{v_0, v_2 \in V_0 \wedge v_1 \in V_1 \wedge \langle v_0, v_2\rangle \in A_0}$, $A_4 = \soast{\langle\langle v_0, v_1\rangle, \langle v_2, v_3\rangle}{v_0, v_2 \in V_0 \wedge v_1, v_3 \in V_1 \wedge \langle v_1, v_3\rangle \in A_1}$, $A_2 = A_3 \cup A_4$, $\pwi{c}{A_3}=\{1\}$ and $\pwi{c}{A_4}=\{2\}$. Clearly the independent sets in $D_0*D_1$ have size at most $mn$ and the only transitive triples in $D_0*D_1$ contain two \one\ arrows which either originate from the same vertex or point at the same. Both types of triples are disagreeable.
\end{proof}

\section{Some tables}

At the end we take the liberty of supplying a table of some Ramsey numbers $r(\alpha, n)$ known today for countable ordinals $\alpha$ and natural numbers $n$. Following the example of Radziszowski in \cite{994R0} we also supply a complementary table containing references. If a reference is given above another, then the upper one gives the upper bound and the lower one gives the lower bound.

\begin{landscape}
\thispagestyle{empty}
\hspace{-6cc}
\begin{array}[b]{c||c*{7}{|c}}&3&4&5&6&7&8&9&m\\\hline\hline
3&6&9&14&18&23&28&36\\\hline
4&9&18&25&&&&\\\hline
\omega&\omega&\omega&\omega&\omega&\omega&\omega&\omega&\omega\\\hline
\omega + 1&\omega 2 + 1&\omega 3 + 1&\omega 4 + 1&\omega 5 + 1&\omega 6 + 1&\omega 7 + 1&\omega 8 + 1&\omega(m-1) + 1\\\hline
\omega + 2&\omega 2 + 4&\omega 3 + 7&\omega 4 + 11&\omega 5 + 16&\omega 6 + 22&\omega 7 + 29&\omega 8 + 37&\omega(m-1) + \frac{m(m-1)}{2}+1\\\hline
\omega + 3&\omega 2 + 8&\omega 3 + 16&&&&&\\\hline
\omega + n&\omega 2 + r(n, 3) + n - 1&&&&&&\\\hline
\lambda 2&\lambda 4&\lambda 8&\lambda 14&\lambda 28&&&\\\hline
\lambda 3&\lambda 9&&&&&&\\\hline
\omega^2&\omega^2&\omega^2&\omega^2&\omega^2&\omega^2&\omega^2&\omega^2&\omega^2\\\hline
\omega^2 + 1&\omega^2 2 + 1&\omega^2 3 + 1&\omega^2 4 + 1&\omega^2 5 + 1&\omega^2 6 + 1&\omega^2 7 + 1&\omega^2 8 + 1&\omega^2(m-1) + 1\\\hline
\omega^2 + 2&\omega^2 2 + 4&\omega^2 3 + 7&\omega^2 4 + 11\qc&\omega^2 5 + 16&\omega^2 6 + 22&\omega^2 7 + 29&\omega^2 8 + 37&\omega^2(m-1) + \frac{m(m-1)}{2}+1\\\hline
\omega^2 + 3&\omega^2 2 + 8&\omega^2 3 + 16&&&&&\\\hline
\omega^2 + n&\omega^2 2 + r(n, 3) + n - 1&&&&&&\\\hline
\omega^2 + \omega&\omega^2 4 + \omega&&&&&&\\\hline
\omega^2 2&\omega^2 10&&&&&&\\\hline
\omega^3&\omega^4&\omega^4&\omega^5&\omega^5&\omega^5&\omega^5&\omega^6&\omega^{2+\lceil\ld(m)\rceil}\\\hline
\omega^3 + n&\omega^4+r(n, 3)&\omega^4 + \omega^3 + r(n, 4) + n - 1&&&&&\\\hline
\omega^4&\omega^7&\omega^7&\omega^{10}&\omega^{10}&\omega^{10}&\omega^{10}&\\\hline
\omega^{5+n}&\omega^{9+2n}&\omega^{9+2n}&\omega^{13+3n}&\omega^{13+3n}&\omega^{13+3n}&\omega^{13+3n}&\omega^{17+4n}&\omega^{1+(4+n)\lceil\ld(m)\rceil}\\\hline
\omega^\omega&\omega^\omega&\omega^\omega&\omega^\omega&\omega^\omega&\omega^\omega&\omega^\omega&\omega^\omega&\omega^\omega\\\hline
\omega^{\omega^2}&\omega^{\omega^2}&\omega^{\omega^2}&&&&&&\\\hline
\kappa\lambda 2&\kappa\lambda 6&&&&&&&\\\hline
\kappa\lambda 3&\kappa\lambda 15&&&&&&&
\end{array}

\newpage

\begin{tabular}[b]{c||c*{7}{
|c
}}
&3&4&5&6&7&8&9&10+m\\\hline\hline
&&&&&\cite{968GY0}&\cite{992MM0}&\cite{982GR0}\\
\rb{3}&&\rb{\cite{955GG0}}&\rb{\cite{955GG0}}&\rb{\cite{964K0}}&\cite{966K0}&\cite{982GR0}&\cite{966K0}\\\hline
&&&\cite{995MR0}&&&&\\
\rb{4}&\rb{\cite{955GG0}}&\rb{\cite{955GG0}}&\cite{965K0}&&&&\\\hline
$\omega$&\cite{930R0}&\cite{930R0}&\cite{930R0}&\cite{930R0}&\cite{930R0}&\cite{930R0}&\cite{930R0}&\cite{930R0}\\\hline
$\omega + 1$&\cite{969HS0}&\cite{969HS0}&\cite{969HS0}&\cite{969HS0}&\cite{969HS0}&\cite{969HS0}&\cite{969HS0}&\cite{969HS0}\\\hline
$\omega + 2$&\cite{969HS0}&\cite{969HS0}&\cite{969HS0}&\cite{969HS0}&\cite{969HS0}&\cite{969HS0}&\cite{969HS0}&\cite{969HS0}\\\hline
$\omega + 3$&\cite{969HS0}&\cite{969HS0}&&&&&&\\\hline
$\omega + n$&\cite{969HS0}&&&&&&&\\\hline
$\lambda 2$&\cite{956ER0}/\cite{974B0}&\cite{964EM0}/\cite{974B0}&\cite{970RP0}/\cite{974B0}&\cite{970RP0}/\cite{974B0}&&&\\\hline
$\lambda 3$&\cite{974B1}/\cite{974B0}&&&&&&\\\hline
$\omega^2$&\cite{957S0}&\cite{957S0}&\cite{957S0}&\cite{957S0}&\cite{957S0}&\cite{957S0}&\cite{957S0}&\cite{957S0}\\\hline
$\omega ^2 + 1$&\cite{969HS0}&\cite{969HS0}&\cite{969HS0}&\cite{969HS0}&\cite{969HS0}&\cite{969HS0}&\cite{969HS0}&\cite{969HS0}\\\hline
$\omega ^2 + 2$&\cite{969HS0}&\cite{969HS0}&\cite{969HS0}&\cite{969HS0}&\cite{969HS0}&\cite{969HS0}&\cite{969HS0}&\cite{969HS0}\\\hline
$\omega ^2 + 3$&\cite{969HS0}&\cite{969HS0}&&&&&\\\hline
$\omega ^2 + n$&\cite{969HS0}&&&&&&\\\hline
$\omega^2 + \omega$&\cite{969HS2}&&&&&&\\\hline
\raisebox{-0.0ex}{$\omega^2 2$}&\cite{969HS1}&&&&&&\\\hline
$\omega^3$&\cite{957S0}&\cite{969HS2}&\cite{974N0}&\cite{974N0}&\cite{974N0}&\cite{974N0}&\cite{974N0}&\cite{974N0}\\\hline
$\omega^3 + n$&\cite{969M0}&\cite{969HS2}&&&&&\\\hline
$\omega^4$&\cite{975N0}&\cite{975N0}&\cite{011HL0}&\cite{011HL0}&\cite{011HL0}&\cite{011HL0}&\\\hline
$\omega^{5+n}$&\cite{979N0}&\cite{979N0}&\cite{979N0}&\cite{979N0}&\cite{979N0}&\cite{979N0}&\cite{979N0}&\cite{979N0}\\\hline
$\omega^\omega$&\cite{973L0}&\cite{973L0}&\cite{973L0}&\cite{973L0}&\cite{973L0}&\cite{973L0}&\cite{973L0}&\cite{973L0}\\\hline
$\omega^{\omega^2}$&\cite{199DL0}&\cite{199DL0}&&&&&&\\\hline
$\kappa\lambda 2$&\this&&&&&&&\\\hline
$\kappa\lambda 3$&\this&&&&&&&
\end{tabular}
\end{landscape}

Combinations of the methods used to find the numbers above would probably generalize some results. 

\section{Problems, Questions \& Conjectures}

To the author the most interesting open question now is the following:

\begin{question}
What is the order of growth of the function $n \mapsto r(I_n, A_3)$?
\end{question}

The author conjectures that he not being able to present a better upper bound is at most partly due to his personal incompetence. So he states the following conjecture---which to his mind seems very modest indeed.

\begin{conjecture}
There exists an $\varepsilon > 0$ such that $\varepsilon n^{2+\varepsilon} < r(I_n, A_3)$ for all natural $n$.
\end{conjecture}

The author was considerably surprised by the second counterexample showing $r(I_2, A_3) > 9$ since there was no \zero\ arrow in it. Hence the following seems to be a natural question to the author:

\begin{question}
Are sharp lower bounds for $r(I_m, A_n)$ for $m$ and $n$ natural numbers always attainable by counterexamples using only \one\ and \two\ arrows?
\end{question}

The author conjectures the answer to this question to be ``No.''.

The author was not able to come up with a formula for the upper bounds given by the recurrence relations for $r(I_m, A_n)$. Perhaps an expert on recurrence relations can find one.

\begin{problem}
Find a nice formula for a general upper bound for $r(I_m, A_n)$!
\end{problem}

It is fairly obvious that the methods of this paper could be extended to determine Ramsey numbers for even more decomposable ordinals. The related finite combinatorics might turn out to be extremely difficult though. We mention just one example as an open problem.

\begin{problem}
Determine $r(\kappa\lambda\omega 2, 3)$ for $\kappa > \lambda$ both weakly compact!
\end{problem}

The motivation behind this research was mainly an interest in countable ordinals. So although we have not discussed Ramsey numbers involving ordinals in $\omega_1\setminus\omega^3$ here at all there are many open problems there. As a teaser we mention the following:

\begin{problem}
Determine $r(\omega^{\omega 2}, 3)$!
\end{problem}

The author in vain tried to strengthen theorem \ref{MAtheorem} by replacing $\MA_{\aleph_1}$ by ``$\omega_1\omega$ is a partition ordinal.''. So he has to leave this as an open question.

\begin{question}
Does $r(\omega_1\omega m, n) = \omega_1\omega r(I_m, S_n)$ whenever $\omega_1\omega$ is a partition ordinal?
\end{question}

A variant of this problem is the following:

\begin{question}
Does $r(\omega_1\omega m, n) = \omega_1\omega r(I_m, S_n)$ whenever $\omega_1\omega^2$ is a partition ordinal?
\end{question}

At the end we want to mention that there are problems in Ramsey theory which to the author's mind could have been attacked algorithmically but were not. The author conjectures that $r(I_4, L_3)$ and $r(I_3, L_4)$ could easily be algorithmically determined. The same should hold true for $r(I_2, S_4)$. Meanwhile the author is rather sceptical about the possibilities to determine $r(I_3, A_3)$ or $r(I_2, A_4)$ without considerable further mathematical insight.

The author wishes to thank at this point Rene Schipperus for getting him interested in this subject, his PhD advisor Stefan Geschke for allowing him to venture into the realms of it, Benjamin Seyfferth for his help in reading parts of the french articles \cite{969HS0}, \cite{969HS1}\ox\ and \cite{969HS2}, Andr\'{e}s Caicedo, Peter Komjath, Lionel Ngyuen Van The, Philipp Schlicht, Stevo Todor\v{c}evi\'{c}, Wolfgang Wohofsky and in particular Jean Larson for discussing this work with them. Moreover \cite{998CG0} was a great source of inspiration for the author. At last I would also like to thank the referee her/his input.\footnote{The research presented in this article was supported by the DFG project

``Continuous Ramsey theory in higher dimension''.}

\small
\bibliography{thilo}

\end{document}